
\def\LaTeX{%
  \let\Begin\begin
  \let\End\end
  \let\salta\relax
  \let\finqui\relax
  \let\futuro\relax}

\def\UK{\def\our{our}\let\sz s}
\def\USA{\def\our{or}\let\sz z}

\UK 



\LaTeX

\USA


\salta

\documentclass[twoside,12pt]{article}
\setlength{\textheight}{24cm}
\setlength{\textwidth}{16cm}
\setlength{\oddsidemargin}{2mm}
\setlength{\evensidemargin}{2mm}
\setlength{\topmargin}{-15mm}
\parskip2mm


\usepackage[usenames,dvipsnames]{color}
\usepackage{amsmath}
\usepackage{amsthm}
\usepackage{amssymb,bbm}
\usepackage[mathcal]{euscript}

\usepackage{cite}
\usepackage{hyperref}
\usepackage[shortlabels]{enumitem}

\usepackage[ulem=normalem,draft]{changes}
%
%

%
 
\definecolor{viola}{rgb}{0.3,0,0.7}
\definecolor{ciclamino}{rgb}{0.5,0,0.5}
\definecolor{blu}{rgb}{0,0,0.7}
\definecolor{rosso}{rgb}{0.85,0,0}

\def\juerg #1{{\color{blue}#1}}
\def\an #1{{\color{cyan}#1}}

\def\pier #1{{\color{red}#1}} 
\def\last #1{{\color{rosso}#1}}
\def\revis #1{{\color{rosso}#1}}

\def\juerg #1{{#1}}
\def\an #1{{#1}}
\def\pier #1{{#1}}
\def\last #1{{#1}}
\def\revis #1{{#1}}




\bibliographystyle{plain}


%

\finqui

\def\Beq{\Begin{equation}}
\def\Eeq{\End{equation}}

\def\Bthm{\Begin{theorem}}
\def\Ethm{\End{theorem}}
\def\Blem{\Begin{lemma}}
\def\Elem{\End{lemma}}
\def\Bprop{\Begin{proposition}}
\def\Eprop{\End{proposition}}

\def\Brem{\Begin{remark}\rm}
\def\Erem{\End{remark}}

\def\Bdim{\Begin{proof}}
\def\Edim{\End{proof}}
\def\Bcenter{\Begin{center}}
\def\Ecenter{\End{center}}
\let\non\nonumber




\def\step #1 \par{\medskip\noindent{\bf #1.}\quad}
\def\jstep #1: \par {\vspace{2mm}\noindent\underline{\sc #1 :}\par\nobreak\vspace{1mm}\noindent}

\def\aand{\quad\hbox{and}\quad}
\def\Lip{Lip\-schitz}
\def\Holder{H\"older}
\def\Frechet{Fr\'echet}
\def\Poincare{Poincar\'e}
\def\lhs{left-hand side}
\def\rhs{right-hand side}



\def\nbh{neighb\our hood}


\def\multibold #1{\def\arg{#1}%
  \ifx\arg\pto \let\next\relax
  \else
  \def\next{\expandafter
    \def\csname #1#1\endcsname{{\boldsymbol #1}}%
    \multibold}%
  \fi \next}

\def\pto{.}

\def\multical #1{\def\arg{#1}%
  \ifx\arg\pto \let\next\relax
  \else
  \def\next{\expandafter
    \def\csname cal#1\endcsname{{\cal #1}}%
    \multical}%
  \fi \next}

\def\multigrass #1{\def\arg{#1}%
  \ifx\arg\pto \let\next\relax
  \else
  \def\next{\expandafter
    \def\csname grass#1\endcsname{{\mathbb #1}}%
    \multigrass}%
  \fi \next}


\def\multimathop #1 {\def\arg{#1}%
  \ifx\arg\pto \let\next\relax
  \else
  \def\next{\expandafter
    \def\csname #1\endcsname{\mathop{\rm #1}\nolimits}%
    \multimathop}%
  \fi \next}

\multibold
qweryuiopasdfghjklzxcvbnmQWERTYUIOPASDFGHJKLZXCVBNM.  

\multical
QWERTYUIOPASDFGHJKLZXCVBNM.

\multigrass
QWERTYUIOPASDFGHJKLZXCVBNM.

\multimathop
diag dist div dom mean meas sign supp .

\def\Span{\mathop{\rm span}\nolimits}


\def\accorpa #1#2{\eqref{#1}--\eqref{#2}}
\def\Accorpa #1#2 #3 {\gdef #1{\eqref{#2}--\eqref{#3}}%
  \wlog{}\wlog{\string #1 -> #2 - #3}\wlog{}}


\def\separa{\noalign{\allowbreak}}

\def\somma #1#2#3{\sum_{#1=#2}^{#3}}

\def\graffe #1{\mathopen\{#1\mathclose\}}
\def\<#1>{\mathopen\langle #1\mathclose\rangle}
\def\norma #1{\mathopen \| #1\mathclose \|}

\def\aeQ{\checkmmode{a.e.\ in~$Q$}}

\def\aet{\checkmmode{a.e.\ in~$(0,T)$}}
\def\aat{\checkmmode{for a.a.\ $t\in(0,T)$}}

\let\hat\widehat
\def\cpto{\,\cdot\,}

\def\iot {\int_0^t}
\def\ioT {\int_0^T}
\def\intQt{\int_{Q_t}}
\def\intQ{\int_Q}
\def\iO{\int_\Omega}

\def\dt{\partial_t}
\def\dn{\partial_{\nn}}
\def\ddt{\frac d{dt}}

\let\emb\hookrightarrow
\def\cpto{\,\cdot\,}

\def\checkmmode #1{\relax\ifmmode\hbox{#1}\else{#1}\fi}


\let\erre\grassR
\let\enne\grassN
\def\erren{\erre^n}




\def\genspazio #1#2#3#4#5{#1^{#2}(#5,#4;#3)}
\def\spazio #1#2#3{\genspazio {#1}{#2}{#3}T0}

\def\L {\spazio L}
\def\H {\spazio H}

\def\C #1#2{C^{#1}([0,T];#2)}


\def\Lx #1{L^{#1}(\Omega)}
\def\Hx #1{H^{#1}(\Omega)}
\def\Wx #1{W^{#1}(\Omega)}

\def\LQ #1{L^{#1}(Q)}

\def\Luno{\Lx 1}
\def\Ldue{\Lx 2}
\def\Linfty{\Lx\infty}

\def\Huno{\Hx 1}
\def\Hdue{\Hx 2}



\let\eps\varepsilon
\let\badphi\phi
\let\phi\varphi

\let\TeXchi\chi                         
\newbox\chibox
\setbox0 \hbox{\mathsurround0pt $\TeXchi$}
\setbox\chibox \hbox{\raise\dp0 \box 0 }
\def\chi{\copy\chibox}



\def\Beta{\widehat\beta}

\def\ej{e_j}
\def\ei{e_i}
\def\lambdaj{\lambda_j}
\def\Vn{V_n}

\def\phin{\phi_n}
\def\mun{\mu_n}
\def\wn{w_n}
\def\phinj{\phi_{nj}}
\def\munj{\mu_{nj}}
\def\un{u_n}
\def\vn{v_n}
\def\zn{z_n}

\def\aj{a_j}
\def\bj{b_j}
\def\cj{c_j}

\def\Tn{T_n}
\def\bphin{\boldsymbol{\phi}_n}
\def\bmun{\boldsymbol{\mu}_n}
\def\bun{\boldsymbol u_n}

\def\Pn{\grassP_n}

\def\ubar{\overline u}
\def\vbar{\overline v}
\def\zbar{\overline z}
\def\hbar{\overline h}
\def\zetabar{\overline\zeta}
\def\phibar{\overline\phi}

\def\phizbar{\overline\phiz}
\def\phinbar{\overline\phin}
\def\munbar{\overline\mun}
\def\psibar{\overline\psi}

\def\ustar{u^*}
\def\phistar{\phi^*}
\def\mustar{\mu^*}
\def\wstar{w^*}
\def\xistar{\xi^*}
\def\zetastar{\zeta^*}

\def\phiz{\phi_0}

\def\Vp{{V^*}}
\def\Wp{W^*}

\def\soluz{(\phi,\mu,w)}
\def\soluzstar{(\phistar,\mustar,\wstar)}
\def\soluzn{(\phin,\mun)}
\def\soluzl{(\psi,\eta,\omega)}
\def\soluza{(p,q,r)}

\def\phih{\phi^h}
\def\muh{\mu^h}
\def\wh{w^h}
\def\psih{\psi^h}
\def\etah{\eta^h}
\def\omegah{\omega^h}

\def\phiQ{\badphi_Q}
\def\phiO{\badphi_\Omega}
\def\wQ{w_Q}
\def\Uad{\calU_{\rm ad}}
\def\UR{\calU_R}
\def\umin{u_{\rm min}}
\def\umax{u_{\rm max}}

\def\normaV #1{\norma{#1}_V}
\def\normaW #1{\norma{#1}_W}
\def\normaVp #1{\norma{#1}_*}

\def\CO{C_\Omega}
\def\cdelta{c_\delta}

\def\Cstar{\an{C^*}}


\usepackage{amsmath}
\DeclareFontFamily{U}{mathc}{}
\DeclareFontShape{U}{mathc}{m}{it}%
{<->s*[1.03] mathc10}{}

\DeclareMathAlphabet{\mathscr}{U}{mathc}{m}{it}

\Begin{document}


%
\title{
Curvature effects in pattern formation: 
\\
well-posedness and optimal control of  a sixth-order Cahn--Hilliard equation}

\author{}
\date{}
\maketitle
\Bcenter
\vskip-1.5cm
{\large\sc Pierluigi Colli$^{(1)}$}\\
{\normalsize e-mail: {\tt pierluigi.colli@unipv.it}}\\[0.25cm]
{\large\sc Gianni Gilardi$^{(1)}$}\\
{\normalsize e-mail: {\tt gianni.gilardi@unipv.it}}\\[0.25cm]
{\large\sc Andrea Signori$^{(2)}$}\\
{\normalsize e-mail: {\tt andrea.signori@polimi.it}}\\[0.25cm]
{\large\sc J\"urgen Sprekels$^{(3)}$}\\
{\normalsize e-mail: {\tt juergen.sprekels@wias-berlin.de}}\\[.5cm]
$^{(1)}$
{\small Dipartimento di Matematica ``F. Casorati'', Universit\`a di Pavia}\\
{\small and Research Associate at the IMATI -- C.N.R. Pavia}\\
{\small via Ferrata 5, I-27100 Pavia, Italy}\\[.3cm] 
$^{(2)}$
{\small Dipartimento di Matematica, Politecnico di Milano}\\
{\small via E. Bonardi 9, I-20133 Milano, Italy}
\\[.3cm] 
$^{(3)}$
{\small Department of Mathematics}\\
{\small Humboldt-Universit\"at zu Berlin}\\
{\small Unter den Linden 6, D-10099 Berlin, Germany}\\
{\small and}\\
{\small Weierstrass Institute for Applied Analysis and Stochastics}\\
{\small Mohrenstrasse 39, D-10117 Berlin, Germany}
\\[10mm]

\Ecenter
\Begin{abstract}
\noindent
This work investigates the well-posedness and optimal control of a \juerg{sixth}-order Cahn--Hilliard equation, 
a higher-order variant of the celebrated and well-established Cahn--Hilliard equation. 
The equation is  endowed with a  source term, where the control variable enters as a distributed mass regulator.
The inclusion of additional spatial derivatives in the \juerg{sixth}-order formulation enables 
the model to capture curvature effects, leading to a more accurate depiction of \juerg{isothermal} phase 
separation dynamics in complex materials systems.
We provide a well-posedness result for the aforementioned system when the corresponding nonlinearity of double-well 
shape is regular and then analyze a corresponding optimal control problem. For the latter, existence of optimal controls is 
established, and the first-order \juerg{necessary} optimality conditions are characterized via a suitable variational inequality. 
These \juerg{results} aim at contributing to improve the  understanding of the mathematical properties and control aspects of the 
\juerg{sixth}-order Cahn--Hilliard equation, offering potential applications in the design and optimization of 
materials with tailored microstructures and properties.

\vskip3mm
\noindent {\bf Key words:} 
{\juerg{Sixth}-order Cahn--Hilliard equation, functionalized Cahn--Hilliard equation, Willmore regularization, curvature effects,
well-posedness, optimal control\juerg{, first-order necessary optimality conditions}.}

\vskip3mm
\noindent 
{\bf AMS (MOS) Subject Classification:} {
		35K55, 
        35K51, 
		49J20, 
		49K20, 
		49J50.  
		}
\End{abstract}

\pagestyle{myheadings}
\newcommand\testopari{\sc Colli -- Gilardi -- Signori -- Sprekels}
\newcommand\testodispari{\sc \an{Analysis and control of a sixth-order Cahn--Hilliard type system}}
\markboth{\testopari}{\testodispari}
%

\section{Introduction}
\label{INTRO}
\setcounter{equation}{0}

The sixth-order Cahn--Hilliard equation represents an extension of the classical Cahn-Hilliard equation that accounts for
curvature effects and higher-order variations of the order parameter. Besides, it provides a more accurate description of
complex \juerg{materials} systems with intricate interface structures. By employing \juerg{a} rigorous analysis, our aim is 
to contribute to the comprehension of the evolution of complex materials, thereby paving the way for advanced applications 
in \juerg{materials} science.

Let $\Omega \subset \mathbb{R}^d$, $d \in \{2,3\}$, be a bounded and smooth domain and let $T>0$. 
The initial-boundary value problem \juerg{under investigation reads~as follows:}
\begin{alignat}{2}
  & \dt\phi - \div(m(\phi)\nabla\mu )
  = S(\phi,u)
  := u - \sigma\phi 
  \qquad && \text{in $Q$,}
  \label{sys:abs:1}
  \\
  & \mu
  = - \eps \Delta w + \tfrac 1\eps f'(\phi) w + \nu w
  \quad && \text{in $Q$,}
  \label{sys:abs:2}
  \\
  & w = -  \revis{\eps \Delta\phi +  \tfrac 1\eps f(\phi)}
  \quad && \text{in $Q$,}
  \label{sys:abs:3}
  \\
  & \dn\phi
  = \pier{ m(\phi)\nabla\mu \cdot \nn}
  = \dn w
  = 0 
  \quad && \text{on $\Sigma$,}
  \label{sys:abs:4}
  \\
  & \phi(0) = \phiz
  \quad && \text{in $\Omega$}.
  \label{sys:abs:5}
\end{alignat}
\Accorpa\Sysabs  {sys:abs:1} {sys:abs:5}
Here, the parabolic cylinder  $Q$ and its boundary $\Sigma$ are defined~by
\Beq
  Q := \Omega\times(0,T)
  \aand
  \Sigma := \Gamma\times(0,T),
  \label{defQS}
\Eeq
where $\Gamma:=\partial\Omega$ denotes the boundary of $\Omega$ that is associated with a \pier{unit} normal vector field $\nn$ and \juerg{the}
outward normal derivative $\dn$.
The above system \Sysabs\ can be viewed as a variant of the classical fourth-order Cahn--Hilliard equation \cite{CH} in the unknowns $\phi$,
 $\mu$, and $w$. The original Cahn--Hilliard equation is widely used to describe \juerg{isothermal} separation processes in binary
 mixtures. The variable $\phi$ represents the local proportion of one of the two components in the binary material and serves as an 
order parameter. To simplify the analysis, it is usually normalized in such a way that the pure states correspond to $\phi = \pm 1$, while 
$\juerg{\{ -1 < \phi < 1\}} $ represents the diffuse interface that occurs in a $\eps$-tubular neighborhood of the interface, with \last{thickness parameter} $\eps>0$. The variable $\mu$ in equation \eqref{sys:abs:2} is referred to as the {\it chemical potential} and corresponds to the first variation of the free energy $\cal E$. Similarly, $w$ is the first variational 
derivative of the Ginzburg--Landau free energy. Namely, it holds that $\mu= \frac {\delta {\cal E}}{\delta \phi}$ and $w=\frac {\delta {\cal G}}{\delta \phi}$, where
\begin{align}
 {\cal E}(\phi) 
 & := 
 {\cal F}(\phi) + \nu {\cal G}(\phi)
 =
\frac 12 \iO \revis{{}\Big(\!-\eps \Delta\phi + \frac 1\eps f(\phi)\Big)^2{}} + \nu \iO \, \Big(\frac \eps2|\nabla\phi|^2 +\frac 1\eps F(\phi)\Big)\,,
  \label{defE}
\end{align}
with natural definition of $\cal F$ and $\cal G$.
It is worth pointing out that $\cal E$ is a higher-order extension of the Ginzburg--Landau free energy $\cal G$ which is associated 
\juerg{with the classical} Cahn--Hilliard equation.
In \eqref{defE}, $F$ indicates \juerg{a double-well shaped nonlinear potential}, and $f$ indicates its derivative. 
A~prototype of $F$ is the {\it classical regular potential} given~by
\Beq
  F(s) := \frac 14 (s^2-1)^2
  \quad \hbox{for \last{every} $s\in\erre$}.
  \label{regpot}
\Eeq
Next, the function $u$ appearing in~\pier{\eqref{sys:abs:1}} is a prescribed distributed function that will play the role of 
\juerg{the} control in the second part of our investigation.
Finally, $\sigma$ and $\nu$ are real constants, with $\sigma$ positive, and $\phiz$~is a prescribed initial datum.

From \eqref{sys:abs:1}, we realize that the mass flux is assumed
to be proportional to the gradient of the chemical potential $\mu$ through the mobility function $m$.
This relationship leads to the variational structure:
\begin{align*}
	\dt \phi = \div (m(\phi) \nabla \mu)+ S(\phi,u),
	\quad 
	\text{with}
	\quad
	\mu= \frac {\delta {\cal E}}{\delta \phi}.
\end{align*}
As a consequence
of the no-flux boundary condition in \eqref{sys:abs:4}, \juerg{integration of} the equation over $\Omega$ leads to  the ODE relation
\begin{align*}
	\frac d {dt} \Big(\, \frac 1 {|\Omega|}\iO \phi(t)\Big)= \frac 1 {|\Omega|} 
	\iO S(\phi(t), u(t)) \quad \text{for all $t \in [0,T],$}
\end{align*}
which \juerg{shows} that the mass dynamics of the order parameter $\phi$  is ruled by the source term $S(\phi, u)$.
 This represents a novelty \juerg{in comparison with} previous works on the system, where $S\equiv 0$ 
\juerg{which results} in mass conservation \last{(see, e.g., \cite{SW, M1, M2})}.

In the Cahn--Hilliard context, the free energy reduces to $\cal G$ with $\nu=1$. On the other hand, considering higher-order terms like the ones in $\cal F$, it \juerg{also can make} sense to allow $\nu$ to be \juerg{negative}.
It turns out that the value and sign of $\nu$ in  \eqref{defE} are essential considerations in both modeling and practical applications involving deformations of elastic vesicles under volume and surface constraints. Specifically, when $\nu=0$, the energy function 
$\mathcal{E}$ simplifies to the well-established \juerg{{\em Willmore functional}} within the phase-field formulation\pier{\cite{CL1,CL2}}. This 
simplification effectively captures the Canham--Helfrich bending energy of surfaces, as demonstrated in, e.g., \cite{DLW,DLRW}.
Besides, when $\nu>0$, the energy function $\mathcal{E}$ is associated with the Willmore regularization of the Ginzburg--Landau energy $\mathcal{G}$. This regularization was employed, for instance, in \cite{TLVW, BCMS}, to investigate anisotropy effects arising during the growth and coarsening of thin films.
Finally, the energy known as the \juerg{\em functionalized Cahn--Hilliard} (FCH) free energy arises when the parameter $\nu$ takes negative values.
This formulation of the energy was primarily developed for mixtures characterized by an amphiphilic structure and nanoscale variations in functionalized polymer chains: see, e.g., \cite{GS, PB}. 
Extensive research has been conducted on the FCH energy, covering topics such as minimization problems, 
bilayer structures, pearled patterns, and network bifurcations. Without the claim of being exhaustive, we refer to \cite{PQ, DP, PQ2, 
DHPW} and the references therein.  
The study of \juerg{sixth}-order Cahn--Hilliard equations has sparked significant interest, and the applications are numerous:
in relation to the dynamics of oil-water-surfactant mixtures, we refer to \cite{PZ1, PZ2, SP}; the faceting of growing surfaces is explored
in works such as \cite{KNR, KR}, while the phase-field-crystal equation is investigated in references like \cite{GW1, GW2, M1, M2, WW}.
\juerg{The novelty in} our contribution is the introduction of a source term $S$ into the system \juerg{that accounts} for mass transfer and exchange \juerg{and leads to a deviation} from the standard assumption of mass conservation. In addition, \juerg{we introduce a control variable in the definition of the source} whose inclusion allows us to manipulate and regulate the behavior of the source term.

\pier{Since} we are not interested in the asymptotic behavior as the interface thickness parameter approaches zero, we set for convenience $\eps=1$.
Moreover, due to our interest in the optimal control application, we consider from the very beginning \juerg{a} constant mobility function 
$m$, that is, we set $m$ to unity for simplicity. Thus, the overall system we are going to analyze is the following:
\begin{alignat}{2}
  & \dt\phi - \Delta\mu 
  = S(\phi,u)
  := u - \sigma\phi 
  \qquad && \text{in $Q$,}
  \label{Iprima}
  \\
  & \mu
  = - \Delta w + f'(\phi) w + \nu w
  \quad && \text{in $Q$,}
  \label{Iseconda}
  \\
  & w = - \Delta\phi + f(\phi)
  \quad && \text{in $Q$,}
  \label{Iterza}
  \\
  & \dn\phi
  = \dn\mu
  = \dn w
  = 0 
  \quad && \text{on $\Sigma$,}
  \label{Ibc}
  \\
  & \phi(0) = \phiz
  \quad && \text{in $\Omega$.}
  \label{Icauchy}
\end{alignat}
\Accorpa\Ipbl Iprima Icauchy
As anticipated above, \juerg{we address in the subsequent stage of our study} an associated distributed control problem, 
where $u$ acts as the control variable.
The \juerg{tracking}-type {\it cost functional} under consideration is given by
\begin{align}
	\calJ (u,\phi)
	 & := 
	\frac {\alpha_1}2 \intQ |\phi - \phiQ|^2
	+ \frac {\alpha_2}2  \iO |\phi(T) - \phiO|^2
	+ \frac {\alpha_3}2 \intQ |u|^2\,,
	\label{Icost}
\end{align}
where the coefficients ${\alpha_i}$, $i=1,2,3$, are nonnegative real numbers 
(not all zero, to avoid a trivial situation),
and $\phiQ$ and $\phiO$ are given target functions {defined in $Q$ and $\Omega$, respectively.}
The set of admissible controls {is} given~by
\Beq
  \Uad := \graffe{u\in\LQ\infty:\ \umin\leq u\leq\umax\ \aeQ},
  \label{IUad}
\Eeq
where $\umin,\umax\in\LQ\infty$ are prescribed functions satisfying $\umin\leq\umax$ \aeQ.
\juerg{The control problem then} consists in minimizing $\calJ$ on $\Uad$
under the \juerg{constraint} that $\phi$ is the first component of the solution $\soluz$ to the state system associated with~$u$.
To summarize, \juerg{we study the following optimal control problem}: 
\begin{align*}
	\text{\bf (P)}
	\quad 
	\min_{u \in \Uad} {\calJ} (u,\phi)
	\quad 
	\text{subject to the constraint that $(\phi,\mu,w)$ solves \Ipbl.}
\end{align*}

The plan of the paper is as follows.
In the upcoming section, we provide a comprehensive list of the specific assumptions we make and present our results.
The well-posedness of the system \Ipbl\ is then studied in Section~\ref{WELLPOSEDNESS}, 
and the control problem is investigated in the last Section~\ref{CONTROL}.


\section{Notation, assumptions and results}
\label{STATEMENT}
\setcounter{equation}{0}

\an{To begin with, let us introduce some notation and conventions. Let
$\Omega$ be} an open set in $\erre^d$, with \an{either} $d=2$ or $d=3$,
which is assumed to be bounded, connected and smooth, \juerg{and whose} Lebesgue measure is denoted by~$|\Omega|$.
\juerg{We denote the outward unit normal field on the boundary $\Gamma:=\partial\Omega$ by $\nn$ and the corresponding 
outward normal derivative
by $\dn$}.
We fix a final time $T>0$, recall the definitions~\eqref{defQS} of $Q$ and $\Sigma$\an{,} and~set
\Beq
	Q_t := \Omega \times (0,t)
	\quad \hbox{for $t\in(0,T]$} \,.
	\label{defQt}
\Eeq
In this paper, \juerg{we employ for any Banach space $X$} the notation $\norma\cpto_X$, $X^*$, and $\< \cpto,\cpto >_X$, to indicate 
the corresponding norm, its dual space, and the related duality pairing between $X^*$ and~$X$.
The only exceptions \juerg{from} this notation for norms \juerg{are} given by the space $H$ introduced below
and the standard Lebesgue spaces $\Lx p$ with $p\in[1,+\infty]$,
whose norms are denoted by $\norma\cpto$ and $\norma\cpto_p$, respectively.
Moreover, the symbol $\norma\cpto_\infty$ might denote the norm in $\LQ\infty$ as well, if no confusion \last{may} arise.
For simplicity, we use the same symbol for the norm \an{in} some space and that in any power \an{thereof}.
Accordingly, when dealing with vector\an{-}valued functions like the gradient of some \an{scalar} function,
we adopt \an{shorthands} like $\L pX$ or $\H1X$ to denote \an{Bochner spaces involving powers of $X$}.
Next, besides the space $H$ announced before, we introduce two further spaces $V$ and~$W$.
Indeed, we~set
\Beq
  H := \Ldue , \quad  
  V := \Huno,
  \aand
  W := \graffe{v\in\Hdue: \ \dn v=0 \hbox{ on $\,\Gamma$}} .
 \label{defspazi}
\Eeq
\juerg{Norm and inner product are in the special case $H$ indicated by $\|\,\cdot\,\|$ and $(\,\cdot\,,\,\cdot\,)$, respectively.}
In connection with the above spaces, we adopt the usual framework of Hilbert triplets \juerg{by identifying $H$ and $\Vp$ with}
subsets of $\Vp$ and $\Wp$ \an{in the usual way.
Namely, we have that}
\begin{align}
  & \< z,v >_V = \textstyle\iO zv
  \aand
  \< z,v >_W = \< z,v >_V,
  \non
  \\
  & \quad \hbox{for every $z\in H$, and $v\in V$ and every $z\in\Vp$ and $v\in W$, respectively}.
  \non
\end{align}
\juerg{Then\pier{, we have that} 
\,$W \emb V \emb H \emb \Vp \emb \Wp$, with dense and compact embeddings}. 

Now, we list our assumptions on the structure of the system to be analyzed.
We \juerg{generally assume:}
\begin{align}
  & 
  \last{\sigma \,\in (0,+\infty) ,
  \quad 
  \lambda \,\in [0,+\infty) }
  \aand
  \nu \in \erre.
  \label{hpconstants}
  \\
  & F \in C^4(\erre) \ \hbox{\juerg{can be written as}} \quad
  F(s) = \Beta(s) - \frac \lambda 2 \, s^2\an{,
  \quad s \in \erre}, 
  \quad \hbox{with $\Beta$ convex} .
  \label{hpBeta}
\end{align}
We set
\Beq
  f := F' , \quad \beta := {\Beta}' \aand \gamma := \beta\beta'\an{,}
  \label{defbeta}
\Eeq
\juerg{so} that 
\begin{align*}
	f(s)= F'(s)= \beta (s) - \lambda s
	\aand 
	f'(s)= \beta'(s) - \lambda,
  \quad \hbox{for every $s\in\erre$,}
\end{align*}
and require that
\begin{align}
  & \beta(0) = \beta''(0) = 0,
  \aand
  \beta'''(s) \geq 0 
  \quad \hbox{for every $s\in\erre$,}
  \label{hpbeta}
  \\
  & \lim_{|s|\to+\infty} \frac {\beta'(s)} {|s|} 
  = + \infty,
  \label{hpbetaprimo}
  \\[2mm]
  & |\beta''(s)|
  \leq C_\beta ( |\beta'(s)| + 1 )
  \quad \hbox{for some $C_\beta>0$ and every $s\in\erre$}.
  \label{hpbetasecondo}
\end{align}
\Accorpa\HPstruttura hpconstants hpbetasecondo

\Brem
\label{Remstructure}
\an{Let us remark that our} structural assumptions imply that 
\Beq
  \hbox{$\beta$, $\beta''$ and $\gamma$ are monotone}.
  \label{betagamma}
\Eeq
This is clear for $\beta$ and $\beta''$, since $\Beta$ is convex and $\beta'''\an{{}=\hat \beta^{(4)}}$ is nonnegative.
As for~$\gamma$, notice that $\gamma'$ is nonnegative,
since $\beta$ and $\beta''$ \an{possess} the same sign.
Moreover, we notice that $|\beta(s)|$ must tend to infinity as $|s|$ \an{approaches} infinity \an{as a consequence of \eqref{hpbetaprimo}},
\juerg{that} the same is true for $\Beta$ since it is convex, and \juerg{that} $\beta(0)=0$ implies that $\Beta$ has a minimum point 
at \an{zero}.
Hence, by \juerg{l'H\^opital's} rule, we deduce~that
\Beq
  \lim_{|s|\to\infty} \frac s{\beta(s)}
  = \lim_{|s|\to\infty} \frac 1{\beta'(s)}
  = 0 
  \aand
  \lim_{|s|\to\infty} \frac {s^3}{\Beta(s)}
  = \lim_{|s|\to\infty} \frac {\an{3}s^2}{\beta(s)}
  = \lim_{|s|\to\infty} \frac {\an{6}s}{\beta'(s)}
  = 0 \an{.}
  \non
\Eeq
\an{Thus, two} consequences \an{follow: first, it holds that} $\Beta$ grows faster than $|s|^3$ as $|s|$ tends to infinity,
so that 
\Beq
  \hbox{$F$ is bounded from below}.
  \label{bddbelow}
\Eeq
The second consequence regards the function $g$ defined~by
\Beq
  g(s) := - \lambda s \beta'(s) + (\nu-\lambda) \beta(s) + (\lambda^2-\lambda\nu) s
  \quad \hbox{for \last{every}  $s\in\erre$,}
  \label{defg}
\Eeq
which we have to consider later~on\juerg{. We have that}
\Beq
  \lim_{|s|\to\infty} \frac {g(s)} {\gamma(s)} = 0
  \aand
  \lim_{|s|\to\infty} \frac {sg(s)} {\gamma(s)} = 0 \,.
  \label{gammag}
\Eeq
Finally, we observe that \an{the} above assumptions \an{are not too restrictive and} satisfied by a wide class of smooth potentials
with either polynomial or exponential growth; in particular, \juerg{they} \an{are met} by the classical regular potential given by 
the \an{expression}~\eqref{regpot}.
\Erem

At this point, we are ready to give a precise formulation of the problem \an{\Ipbl\ presented} in the Introduction.
\an{As anticipated, the control variable $u$ plays the role of a prescribed and bounded forcing term in the well-posedness part.}
Given a constant $M>0$, as well as data $\phiz$ and $u$ satisfying
\Beq
  \phiz \in W
  \aand
  u \in\LQ\infty
  \quad \hbox{with} \quad
  \norma u_\infty \leq M\,,
  \label{hpdati}
\Eeq
we look for a triplet $\soluz$ with the regularity 
\begin{align}
  & \phi \in \H1\Vp \cap \L\infty W,
  \label{regphi}
  \\
  & \mu \in \L2V,
  \label{regmu}
  \\
  & w \in \L2W,
  \label{regw}
\end{align}
\Accorpa\Regsoluz regphi regw
that solves the following problem:
\begin{align}
  & \< \dt\phi , v >_\an{V}
  \an{{}+ \iO \nabla\mu \cdot \nabla v
  + \sigma \iO \phi v
  }
  = \iO u v
  \quad \hbox{\aet, for every $v\in V$,}
  \label{prima}
  \\
  & - \Delta w + \beta'(\phi) w + (\nu-\lambda) w
  = \mu
  \quad \aeQ,
  \label{seconda}
  \\
  & - \Delta\phi + \beta(\phi) - \lambda\phi
  = w
  \quad \aeQ,
  \label{terza}
  \\
  & \phi(0) = \phiz \,.
  \label{cauchy}
\end{align}
\Accorpa\Pbl prima cauchy
Notice that the boundary conditions $\dn\phi=0$ and $\dn w=0$ are contained in \eqref{regphi} and \eqref{regw},
while the analogue for~$\mu$ just holds in a generalized sense as a consequence
\pier{of~\eqref{prima}.}

As for the assumption on~$u$,
we observe that $M$ does not play any role if just well-posedness is considered.
However, when dealing with the control problem, we let $u$ vary
and \an{need} uniform bounds for the corresponding solutions,
so that the introduction of $M$ is useful:
\an{it is worth noticing that}
the bounds we are going to find are independent of~$u$ and depend just on~$M$.

The first result of ours deals with well-posedness, regularity, and stability \an{of \juerg{the} system  \Pbl}.

\Bthm
\label{Wellposedness}
Let the assumptions \HPstruttura\ and \eqref{hpdati} on the structure of the system and the data be satisfied.
Then, problem \Pbl\ has a unique solution \an{$(\phi,\mu,w)$} satisfying \Regsoluz.
Moreover, the estimate
\Beq
  \norma\phi_{\H1\Vp\cap\L\infty W}
  + \norma\mu_{\L2V}
  + \norma w_{\L2W}
  \leq K_1 
  \label{stability}
\Eeq
holds true with a positive constant $K_1$ that depends only on $\Omega$, $T$, the structure of the system,
the initial datum $\phiz$ and~$M$.
\Ethm

\Brem
\label{Piureg}
\an{The unique} solution $\soluz$ to problem \Pbl\ satisfying \Regsoluz\ \an{actually} enjoys the further regularity properties
\Beq
  \phi \in \L2{\Hx4}
  \aand
  w \in \L\infty H \cap \L2{\Hx3} ,
  \label{piureg}
\Eeq
as we immediately prove \juerg{by repeated use of} elliptic regularity theory.
Since $\phi$ is \an{uniformly} bounded as a consequence of~\eqref{regphi} and the continuity of the embedding $W\emb\Linfty$,
we see that \eqref{regphi} and \eqref{terza} imply that $w\in\L\infty H$.
Now, we observe that \eqref{regphi} also implies that $|\nabla\phi|^2\in\L\infty H$,
so that $\Delta\beta(\phi)\in\L\infty H$ and $\beta(\phi)\in\L\infty\Hdue$.
Thus, \eqref{terza} and \eqref{regw} yield that $\phi\in\L2{\Hx4}$.
Next, the regularity just obtained, and that of $w$ \an{in \eqref{regw}}, imply that 
\pier{$\nabla\phi\in\L2\Linfty$  
and $\nabla w\in\L2H$}.
Hence, \pier{as $w\in \L\infty H$\last{,} it turns out that} $\nabla(\beta'(\phi)w)\in\L2H$, i.e., $\beta'(\phi)w\in\L2V$.
Therefore, by also accounting for \eqref{regmu}, we conclude from \eqref{seconda} that $w\in\L2{\Hx3}$.
We could continue in improving the regularity of~$\phi$.
However, we just observe one more property of $\nabla\phi$ \juerg{that will be needed} in the sequel.
Namely, \juerg{it holds that}
\Beq
  |\nabla\phi| \in  \L4\Linfty,
  \label{piuregbis}
\Eeq
\juerg{by virtue of the embedding} $\,\L\infty\Huno\cap\L2{\Hx3}\emb\L4\Linfty$.
We conclude this remark by observing that if \eqref{stability} holds true for $\soluz$,
then analogous estimates associated to the regularity~\eqref{piureg} and \eqref{piuregbis} are satisfied as well \an{by the solution},
since each of the above steps also provides a corresponding estimate.
\Erem

The second one regards continuous dependence of the solution \an{with respect to} the control variable~$u$.
\Bthm
\label{Contdep}
Let the assumptions \HPstruttura\ on the structure of the system be fulfilled,
as well as \eqref{hpdati} for $\phiz$,
and let $u_i\in\LQ\infty$, $i=1,2$, \juerg{satisfy} $\norma{u_i}_\infty\leq M$.
If $(\phi_i,\mu_i,w_i)$ are the corresponding solutions \an{to system \Pbl},
\juerg{then} the estimate
\begin{align}
  & \norma{\phi_1-\phi_2}_{\C0V\cap\L2{\Hx4}}
  + \norma{\mu_1-\mu_2}_{\L2H}
  \non
  \\
  & \quad {}
  + \norma{w_1-w_2}_{\L2W}
  \leq K_2 \, \norma{u_1-u_2}_{\L2\Vp}
  \label{contdep}
\end{align}
holds true with a positive constant $K_2$ that depends only on $\Omega$, $T$, the structure of the system,
the initial datum $\phiz$ and~$M$.
\Ethm

As announced in the Introduction, we address an associated distributed control problem \an{for system \Ipbl}.
Here, we make our precise \juerg{assumptions:} 
\begin{align}
  & \hbox{$\alpha_i$ \ are nonnegative real numbers for \ $\an{i=1,2,3.}$}
  \label{hpalpha}
  \\
  & \phiQ \in \LQ2
  \aand
  \phiO \in \an{V}. 
  \label{hpcost}
  \\
  & \umin ,\, \umax \in \LQ\infty
  \quad \hbox{with} \quad
  \umin \leq \umax \quad \aeQ.
  \label{hpUad}
\end{align}
\juerg{Moreover,} we define the cost functional, the control space, and the set of admissible controls, by setting
\begin{align}
	& \calJ (u,\phi)
    := \frac {\alpha_1}2 \intQ |\phi - \phiQ|^2
	+ \frac {\alpha_2}2  \iO |\phi(T) - \phiO|^2
	+ \frac {\alpha_3}2 \intQ |u|^2\,,
	\label{cost}
	\\
    & \rlap{$\calU := \LQ\infty,$\aand $\Uad := \graffe{\an{u \in {\cal U}\,:}\, \umin\leq u\leq\umax\ \aeQ}$.}
    \label{Uad}
\end{align}
\an{Of course, for \eqref{cost} to be meaningful, it  \juerg{suffices to require} $\phiO \in H$, but that would not be enough to handle the control problem for technical reasons that will be clarified later on (cf. Theorem \ref{Wellposednessa}).}
Thus, the control problem \juerg{under investigation} is the following:
\begin{align}
  & \hbox{Minimize $\calJ(u,\phi)$ under the conditions:}
  \non
  \\
  & \quad \hbox{$u\in\Uad$, and $\phi$ is the first component of the solution $\soluz$}
  \non
  \\
  & \quad \hbox{to the state system \last{\Pbl} associated with~$u$}.
  \label{control}
\end{align}
We have considered \an{a basic} form \an{for} the cost functional \an{for simplicity}.
However, \an{let us claim that we can actually afford to include in the cost functional}
\eqref{cost} \an{an additional term} involving the component $w$ and some target function~$\wQ$
(see the forthcoming Remark~\ref{Moregencost}).

In Section~\ref{CONTROL}\an{,} we prove that \an{there \pier{exist} optimal controls} and \revis{provide first-order} necessary conditions for $\ustar$ to be an optimal control \juerg{(or, more generally, a locally optimal control in the sense of $L^p$ for $1\le p\le +\infty$, 
see the concluding Remark 4.10)}, which \an{reads} as follows:
if $\ustar\in\Uad$ is the optimal control and $\soluzstar$ is the corresponding state,
then the variational inequality
\Beq
  \intQ (\alpha_3 \, \ustar + p) (u-\ustar) \geq 0
  \non
\Eeq
is satisfied for every $u\in\Uad$, where $p$ is the first component of the solution $\soluza$ to a proper adjoint system \an{associated to \Pbl}.
Here, 
\an{for simplicity, we just}
confine ourselves to \an{mention} 
that it is a \an{suitable} weak formulation of the following formal backward\an{-in-time} system
\begin{alignat}{2}
  & - \dt p - \Delta r + \sigma p
  + f''(\phistar) \wstar q - f'(\phistar) r
  = \alpha_1 (\phistar-\phiQ)
  \non
  \quad && \hbox{\an{in $Q$,}}
  \\
  & q = - \Delta p
  \non
  \quad && \hbox{\an{in $Q$,}}
  \\
  & r + \Delta q - \nu q - f'(\phistar) q
  = 0
  \non
  \quad && \hbox{\an{in $Q$,}}
\end{alignat}
complemented \revis{with homogeneous} Neumann boundary conditions for all the \an{variables} and the final condition
\Beq
  p(T) = \alpha_2 (\phistar(T)-\phiO) \,
  \non
   \quad \hbox{\an{in $\Omega$.}}
\Eeq

We conclude this section by collecting some useful tools that will be employed throughout the paper.
Besides the \Holder\ inequality, we often account for the Young, Sobolev and Poincar\'e inequalities 
as well as for some inequalities associated to elliptic regularity theory and to the compact embeddings $V\emb\Lx p$ and $W\emb V$.
\juerg{In fact, we have }
\begin{align}
  & ab \leq \delta a^2 + \frac 1{4\delta} \, b^2
  \quad \hbox{for every $a,b\in\erre$ and $\delta>0$,}
  \label{young}
  \\[2mm]
  \separa
  & \norma v_p \leq \CO \normaV v
  \quad \hbox{for every $v\in V$ and $p\in[1,6],$}
  \label{sobolev}
  \\[2mm]
  \separa
  & \normaV v
  \leq \CO \, \bigl( \norma{\nabla v} + |\vbar| \bigr)
  \quad \hbox{for every $v\in V$,}
  \label{poincare}
  \\[2mm]
  \separa
  & \normaW v \leq \CO \bigl( \norma{\Delta v} + \normaVp v \bigr)
  \label{elliptic}
  \quad \hbox{for every $v\in W$,}
  \\[2mm]
  & \norma v_p
  \leq \delta \, \norma{\nabla v} + C_{\Omega,p,\delta} \, \normaVp v
  \quad \hbox{for every $v\in V$, $p\in[1,6)$ and $\delta>0$,}
  \label{compact}
  \\[2mm]
  \separa
  & \normaV v
  \leq \delta \, \norma{\Delta v} + C_{\Omega,\delta} \, \normaVp v
  \quad \hbox{for every $v\in W$ and $\delta>0$,}
  \label{compactbis}
\end{align}
where $\vbar$ denotes the mean value of $v$ and $\normaVp\cpto$ \an{a} norm in $\Vp$ \juerg{which will be} introduced below 
\juerg{in} \eqref{normaVp}. 
In the above inequalities, \an{the constant} $\CO$~depends only on~$\Omega$\juerg{, while} $C_{\Omega,\delta}$ 
and $C_{\Omega,p,\delta}$ also depend on $\delta$ and~$(p,\delta)$, respectively.

More generally, \juerg{we define the generalized mean value $\vbar$ of a generic element $v\in\Wp$} 
by \an{setting}
\Beq
  \vbar := \frac 1{|\Omega|} \, \< v , 1 >_W\,,
  \label{defmean}
\Eeq 
where $1$ stands for the constant function that takes the value $1$ everywhere in~$\Omega$.
\last{F}or any $s\in\erre$, we still denote by $s$ the corresponding constant functions \an{in} $\Omega$ and~$Q$.
Notice that the above definition \eqref{defmean} is meaningful, since $1$ actually belongs to~$W$,
and that $\vbar$ reduces to the usual mean value \an{when} $v\in H$.
The same notation $\vbar$ is also employed if $v$ is a time-dependent function.

Next, we recall an important tool which is commonly used when working with problems connected to the Cahn--Hilliard equation.
\juerg{To this end, consider} the weak formulation of the Poisson equation $-\Delta z=\zeta$
with homogeneous Neumann boundary conditions. 
Namely, for a given $\zeta\in\Vp$ (which does not necessarily belong to $H$), we consider the problem
of finding
\begin{align}
  z \in V
  \quad \hbox{such that} \quad
  \iO \nabla z \cdot \nabla v
  = \< \zeta , v >_\an{V}
  \quad \hbox{for every $v\in V$}.
  \label{neumann}
\end{align}
	Since $\Omega$ is connected and smooth, it is well known that the above problem admits \pier{solutions} $z$ if and only if $\zeta$ has zero mean value.
Hence, we can introduce the \pier{following} solution operator $\calN$ by setting
\begin{align}
  & \calN: \dom(\calN) := \graffe{\zeta\in\Vp:\ \zetabar=0} \to \graffe{z\in V:\ \zbar=0},
  \quad 
  {\cal N}: \zeta \mapsto z,
  \label{defN}
\end{align}
where $z$  is the unique solution to \eqref{neumann} \pier{coupled with} $\zbar=0$.
It turns out that $\calN$ is an isomorphism between the above spaces, \an{and it} follows that the formula
\Beq
  \normaVp\zeta^2 := \norma{\nabla\calN(\zeta-\zetabar)}^2 + |\zetabar|^2
  \quad \hbox{for every $\zeta\in\Vp$}
  \label{normaVp}
\Eeq
defines a Hilbert norm in $\Vp$ that is equivalent to the standard dual norm of~$\Vp$.
From the above properties, one can obtain the following identities:
\begin{align}
  & \iO \nabla\calN\zeta \cdot \nabla v
  = \< \zeta , v >_\an{V}
  \quad  \hbox{for every $\zeta\in\dom(\calN)$ and $v\in V$\pier{,}}
  \label{dadefN}
  \\
  & \< \zeta , \calN\xi >_\an{V}
  = \< \xi , \calN\zeta >_\an{V}
  \quad \hbox{for every $\zeta,\xi\in\dom(\calN)$\pier{,}}
  \label{simmN}
  \\
  & \< \zeta , \calN\zeta > _\an{V}
  = \iO |\nabla\calN\zeta|^2
  = \normaVp\zeta^2
  \quad \hbox{for every $\zeta\in\dom(\calN)$.}
  \label{danormaVp}
\end{align}
\juerg{Moreover,} \pier{we point out that}
\Beq
  \< \dt\zeta(t) , \calN\zeta(t) >_V
  = \frac 12 \, \ddt \, \normaVp{\zeta(t)}^2
  \quad \aat\,,
  \label{propN} 
\Eeq
which holds true for every $\zeta\in\H1\Vp$ satisfying $\zetabar=0$ \aet.
\Accorpa\PropN defN propN

Finally, without further reference later on, we are going to employ the following convention: 
the small-case symbol $c$ denotes a generic constant
that depends only on the structure of the system, $\Omega$, $T$, the initial datum~$\phiz$,
and the constant $M$ that appears in~\eqref{hpdati}. 
In particular, the values of $c$ \an{are} independent of $u$ and the \an{approximation} parameter $n$ we introduce in Section~\ref{EXISTENCE}.
Notice that the meaning of $c$ may vary from line to line and even within the same line.
In addition, whenever a positive constant $\delta$ enters the computation, the related symbol~$\cdelta$, in place of a general~$c$, 
denotes constants that depend on~$\delta$, in addition.
We use different notations for precise constants we could refer~to, like, e.g., in~\eqref{sobolev}.


\section{The state system}
\label{WELLPOSEDNESS}
\setcounter{equation}{0}

In this section, we prove Theorems~\ref{Wellposedness} and~\ref{Contdep}.
To this end, it is convenient to observe that our problem can be formulated in an alternative \an{equivalent} way.

\Bprop
\label{Pblbis}
If $\soluz$ satisfies \Regsoluz\ and solves \an{\Pbl},
then the pair $(\phi,\mu)$ also solves the variational equation
\begin{align}
  & \iO \Delta\phi \, \Delta v
  - \iO \Delta\beta(\phi) \, v
  - \iO \beta'(\phi) \Delta\phi\, v
  \non
  \\
  & \quad {}
  + \iO \gamma(\phi) \, v
  + (2\lambda-\nu) \iO \Delta\phi \, v
  + \iO g(\phi) \, v
  \non
  \\
  & = \iO \mu v
  \quad \hbox{\aet, for every $v\in W$,}
  \label{secondabis}
\end{align}
where the function $g:\erre\to\erre$ is defined by \eqref{defg}.
Conversely, if $\phi$ and $\mu$ satisfy \an{the regularities in} \accorpa{regphi}{regmu} 
and the pair $(\phi,\mu)$ solves \eqref{secondabis},
then we have the following:\\
\juerg{{\rm (i)}~The} function $w$ given~by
\Beq
  w := -\Delta\phi + \beta(\phi) - \lambda\phi
  \label{defw}
\Eeq
satisfies \eqref{regw} and the estimate
\Beq
  \norma w_{\L2W} \leq C \bigl( \norma\phi_{\L2W} + \norma\mu_{\L2H} + 1 \bigr)\,,
  \label{stimaw}
\Eeq
with a positive constant $C$ that depends only on $\Omega$, $T$ and the structure of the system\juerg{.
\\
{\rm (ii)}~The} triplet $\soluz$ solves \an{\Pbl}.
\Eprop

\Bdim
We observe once and for all that \eqref{regphi}, the continuous embedding $W\emb\Linfty$, and the smoothness of $\beta$,
\juerg{imply} that $\phi$, $\beta(\phi)$, and $\beta'(\phi)$, \an{are bounded in}~$\LQ\infty$.
\last{Besides, throughout} the proof\an{,} we repeatedly owe to elliptic regularity theory.

Assume first that $\soluz$ satisfies \Regsoluz\ and solves the problem \accorpa{prima}{terza}.
\an{Then, we} eliminate $w$ in \eqref{seconda} by means of \eqref{terza} \an{to} obtain~that
\Beq
  -\Delta \bigl( -\Delta\phi + \beta(\phi) - \lambda\phi \bigr)
  + \bigl( \beta'(\phi) - \lambda + \nu \bigr) \bigl( -\Delta\phi + \beta(\phi) - \lambda\phi \bigr)
  = \mu
  \quad \aeQ\an{.}
  \label{quartordine}
\Eeq
Now, we recall that $w\in\L2W$
\an{and the identity}
\Beq
  \Delta\beta(\phi) = \beta''(\phi) |\nabla\phi|^2 + \beta'(\phi) \Delta\phi\,,
  \non 
\Eeq
\an{where} both $|\nabla\phi|^2$ and $\Delta\phi$ belong to $\L\infty H$,
the former since $W\emb\Wx{1,4}$.
Hence, $\Delta\beta(\phi)$ belongs to $\L\infty H$.
Since\, $\dn\beta(\phi)=\beta'(\phi)\dn\phi=0$ \juerg{on $\Sigma$}, we conclude that $\beta(\phi)\in\pier{\L2 W}$,
\an{whence}, by comparison in \eqref{terza}, we also have that $\Delta\phi\in \juerg{L^2(0,T;W)}$.
Thus, we can distribute the first Laplac\an{i}an in \eqref{quartordine} to the single summands.
By recalling \eqref{defg}\an{,} we obtain~that
\Beq
  \Delta^2 \phi 
  - \Delta\beta(\phi)
  - \beta'(\phi) \Delta\phi
  + \gamma(\phi)
  + (2\lambda-\nu) \Delta\phi
  + g(\phi)
  = \mu
  \quad \aeQ \,.
  \label{secondabisPde}
\Eeq
\an{To} derive~\eqref{secondabis}, it suffices to multiply \eqref{secondabisPde} by an arbitrary $v\in W$,
integrate over $\Omega$ and use \an{integration by parts} in the first integral \an{along with the property that} both $\Delta\phi$ and $v$ have a zero normal derivative \juerg{on $\Sigma$}.

Conversely, assume that $\phi$ and $\mu$ \an{fulfill} \accorpa{regphi}{regmu}
and the pair $(\phi,\mu)$ solves \eqref{secondabis}.
Then\an{, since $\dn\beta(\phi)=0$ \juerg{on $\Sigma$}, it holds that}
\Beq
  \iO \Delta\beta(\phi) \, v
  = \iO \beta(\phi) \Delta v
  \quad \hbox{\aet, for every $v\in W$}.
  \non
\Eeq
Hence, by recalling the definitions \eqref{defbeta} and \eqref{defg} of $\gamma$ and~$g$,
we can rewrite \eqref{secondabis} in the~form
\begin{align}
  & \iO \bigl( -\Delta\phi - \beta(\phi) - \lambda\phi \bigr) (-\Delta v)
  + \iO \bigl( \beta'(\phi) - \lambda + \nu \bigr) \bigl( -\Delta\phi + \beta(\phi) - \lambda\phi \bigr) v
  \non
  \\
  & = \iO \mu v
  \quad \hbox{\aet, for every $v\in W$}.
  \non
\end{align}
Therefore, if $w$ is given by~\eqref{defw},
then $w\in\L\infty H$, and the \an{equation} \eqref{terza} is satisfied.
Moreover, since \eqref{secondabis} has the form
\Beq
  \iO \Delta\phi \, \Delta v
  = \iO h v
  \quad \hbox{\aat, for every $v\in W$,}
  \non
\Eeq
with $h\in\L2H$, we deduce that $\Delta^2\phi\in\L2H$ and $\dn\Delta\phi=0$ \juerg{on $\Sigma$}.
Thus, we also have that $w\in\L2{\Hx2}$ and $\dn w=0$ \juerg{on $\Sigma$}, i.e., \eqref{regw}.
Furthermore, it is clear~that
\Beq
  \iO w (-\Delta v)
  + \iO \bigl( \beta'(\phi) - \lambda + \nu \bigr) \, w \, v
  = \iO \mu v
  \quad \hbox{\aet, for every $v\in W$}\an{,}
  \non
\Eeq
\an{from which we obtain} \eqref{seconda} since $w\in\L2W$.
Finally, we have~that
\Beq
  \norma w_{\L2W}
  \leq c \, \bigl(  \norma\mu_{\L2H} + \norma w_{\L2H} \bigr)
  \leq c \, \bigl(  \norma\mu_{\L2H} + \norma\phi_{\L2W} +1 \bigr)\,,
  \non
\Eeq
that is, \eqref{stimaw} with some \an{computable} constant $C$ as in the statement.
\Edim

\Brem
\label{Pblter}
Since both $\phi$ and $w$ belong to $\L2W$, \juerg{we conclude from \eqref{terza} that} 
\Beq
  \Delta\phi \in \L2W \,.
  \label{regDeltaphi}
\Eeq
Hence, the first integral of \eqref{secondabis} \an{can be rewritten as}  $\iO\Delta^2\phi\,v$,
and \eqref{secondabis} itself is equivalent~to
\an{
\begin{align}
  & \Delta^2\phi - \Delta\beta(\phi) - \beta'(\phi) \Delta\phi
  + \gamma(\phi) 
  + (2\lambda-\nu) \Delta\phi 
  + g(\phi)  \non
  = \mu 
  \quad \aeQ \,,
 \non 
\end{align}}
\pier{that is, \last{equation} \eqref{secondabisPde}.}
\Erem

\juerg{In summary}, the problem \Pbl\ \an{may be reformulated as}:
\Beq
  \hbox{\it Find $(\phi,\mu)$ satisfying \accorpa{regphi}{regmu}, \eqref{prima}, \eqref{secondabis} and \eqref{cauchy}.}
  \label{pblbis}
\Eeq
In this new framework, due to \eqref{stimaw},
the stability estimate \eqref{stability} is equivalent 
to the estimate obtained by \eqref{stability} itself by ignoring the part regarding~$w$.


\subsection{Existence}
\label{EXISTENCE}

\an{In this section, we}
prove the existence and stability part of Theorem~\ref{Wellposedness}
by constructing a solution that satisfies the estimate~\eqref{stability}.
By Proposition~\ref{Pblbis}, we can consider the problem  \an{in the form expressed in} \eqref{pblbis}.
To this end, we introduce a \an{corresponding} \juerg{discretization} by means of a Faedo--Galerkin scheme.

\an{In this direction, let} $\graffe{\lambdaj}_{j\geq1}$ and $\graffe{\ej}_{j\geq1}$
be the sequence of the eigenvalues and an orthonormal system of corresponding eigenfunctions
of the Neumann problem for the Laplace equation, i.e., \juerg{we have} \last{that}
\begin{align}
  & 0 = \lambda_1 < \lambda_2 \leq \lambda_3 \leq \dots
  \aand \lim_{j\to\infty} \lambdaj = + \infty\,,
  \label{eigenvalues}
  \\
  & \ej\in V
  \aand
  \iO \nabla\ej \cdot \pier{\nabla v}
  = \lambda_j \iO \ej v
  \quad \hbox{for every $v\in V$ and $j=1,2,\dots$\,,}
  \label{eigenfunctions}
  \\
  & \iO e_i e_j = \delta_{ij}
  \quad \hbox{for $i,j=1,2,\dots$\,,}
  \aand 
  \hbox{$\graffe{e_j}_{j\geq1}$ is a complete system in $H$\,,}
  \label{orthogonality}
\end{align}
\Accorpa\EigenPbl eigenvalues orthogonality
where $\delta_{ij}$ is the Kronecker symbol.
Notice that \juerg{actually $\ej \in W$,} since $\Omega$ is \last{assumed to be} smooth\last{, and}
$e_1$ is a constant \juerg{eigenfunction}.
We~set
\Beq
  \Vn := \Span \graffe{e_1,\dots,e_n},
  \quad \hbox{for $n=1,2,\dots$},
  \label{defVn}
\Eeq
and observe that the union of these spaces is dense in both $V$ and~$H$.
\an{Besides, we point out once \juerg{and} for all that $V_1=\Span \graffe{e_1}$ consists of the space of constant functions\juerg{, which}
 are therefore admissible test functions in the \juerg{Galerkin} scheme.}

\Brem
\label{RemVn}
We also notice that $\Vn\subset W$ for every~$n$ and~that
\Beq
  \calN v \in \Vn
  \quad \hbox{for every $v\in\Vn$ with zero mean value}.
  \label{calNvn}
\Eeq
Indeed, if \pier{$v\in \Vn $ has zero mean value, then it can be represented as} $v={\somma j2n} \cj\ej$ \pier{for suitable coefficients $\cj$, $j=2, \ldots, n$. On the other hand,} the condition $z=\calN v$ means that
\Beq
  z \in \pier{W}, \quad 
  \zbar = 0
  \aand
 - \Delta z = v \,\juerg{\mbox{ in \,$\Omega$}},
  \label{pier1}
\Eeq
\pier{whence from \eqref{eigenfunctions} it easily follows that $z:={\somma j2n} \lambdaj^{-1}\cj\ej \in \Vn$ solves \eqref{pier1}.}
%

Next, let $\Pn:H\to\Vn$ be the \juerg{orthogonal} projection operator \last{onto $V_n$},
and let $Y$ \juerg{be any} of the spaces $H$, $V$ and~$W$.
Then, we have that 
\Beq
  \norma{\Pn v}\an{_Y} \leq \CO \, \norma v\an{_Y}
  \quad \hbox{for every $v\in Y$,}
  \label{stimaPnH}
\Eeq
where the constant~$\CO$ 
(\an{with} $\CO=1$ if $Y=H$)
depends only on~$\Omega$.
Moreover, if $v\in\L2Y$ and $\vn$ is defined by
$\vn(t):=\Pn(v(t))$ \aat,
then
\Beq
  \norma\vn_{\L2Y} \leq \CO \, \norma v_{\L2Y},
  \aand
  \vn \to v 
  \quad \hbox{strongly in $\L2Y$}.
  \label{convPh}
\Eeq
\an{\juerg{These} properties are straightforward if} $Y=H$\last{, whereas for}
the other cases\an{, we refer\last{, e.g.,} to \cite[Rem~3.3]{CGSS1}.}
\Erem

\an{We are now ready to} state the discrete problem\juerg{: find a pair $\soluzn$ 
with} 
\begin{align}
  & \phin \in \H1\Vn
  \aand
  \mun \in \L2\Vn,
  \label{regsoluzn}
  \\[2mm]
  \separa
  & \iO \dt\phin \, v 
  \an{{}+ \iO \nabla\mun \cdot \nabla v
  + \sigma \iO \phin v
  }
  = \iO u v
  \quad \hbox{\aet, for every $v\in\Vn$,}
  \label{priman}
  \\
  \separa
  & \iO \Delta\phin \, \Delta v
  - \iO \Delta\beta(\phin) \, v
  - \iO \beta'(\phin) \Delta\phin \, v
  \non
  \\
  & \quad {}
  + \iO \gamma(\phin) \, v
  + (2\lambda-\nu) \iO \Delta\phin \, v
  + \iO g(\phin) \, v
  \non
  \\
  & = \iO \mun \, v
  \quad \hbox{\aet, for every $v\in\Vn$,}
  \label{secondan}
  \\
  \separa
  & \iO \phin(0) \, v
  = \iO \phiz \, v
  \quad \hbox{for every $v\in\Vn$},
  \label{cauchyn}
\end{align}
\pier{where the last condition implies that $\phin(0)  = \Pn(\phiz)$.}
\Accorpa\Pbln regsoluzn cauchyn

\Bthm
\label{Wellposednessn}
\an{For every $n \in \enne$, t}here exists a unique pair $\soluzn$ satisfying \Pbln.
\Ethm

\Bdim
\an{First of all, let us fix $n \in \enne$.}
We expand the \juerg{sought} solution \juerg{in the form}
\Beq
  \phin(t) = \somma j1n \phinj(t) \ej\,,\quad
  \mun(t) = \somma j1n \munj(t) \ej\,,
  \non
\Eeq
and introduce \an{the real unknowns of the problem, which are} the vector\an{-}valued functions $\bphin$, $\bmun$ and $\bun$ (to be \last{intended} as columns), 
by setting\an{,} \aat\an{,}
\begin{align}
  & \bphin(t) := (\phinj(t))_{j=1}^n  \,, \quad
  \bmun(t) := (\munj(t))_{j=1}^n
  \aand
  \bun(t) := (u_{nj}(t))_{j=1}^n\,,
  \non
  \\
  & \quad \hbox{where} \quad
  u_{nj}(t) = \iO u(t) \ej
  \quad \hbox{for $j=1,\dots,n$}.
  \non
\end{align}
We aim at presenting the above problem as a system of ordinary differential equations in the unknown $(\bphin,\bmun)$.
To this end, we recall \an{the identity}
\Beq
  - \Delta\beta(\phin)
  = - \beta''(\phin) |\nabla\phin|^2 - \beta'(\phin) \Delta\phin\an{.}
  \non
\Eeq
\an{From this,} for $i=1,\dots,n$, we have that
\begin{align}
  & \iO \bigl( - \Delta\beta(\phin) - \beta'(\phin) \Delta\phin \bigr) \, \ei
  \non
  \\
  & = - \iO \beta''(\phin) \Bigl| \somma j1n \phinj \nabla\ej \Bigr|^2 \, \ei
  - 2 \iO \beta'(\phin) \somma j1n \phinj \lambdaj \ej \ei \,.
  \non
\end{align}
Then, the system \accorpa{regsoluzn}{secondan} is equivalent~to \an{the vector-valued system}
\begin{align}
  & \bphin \in \H1\erren
  \aand
  \bmun \in \L2\erren\,,
  \label{regode}
  \\[1mm]  
  & \bphin' \an{{}+ A \bmun + \sigma\bphin } = \bun\,,
  \label{odea}
  \\[1mm]
  & B^2 \bphin + \calB(\bphin) - 2 \calF(\bphin) + \calG(\bphin) - (2\lambda-\nu) B \bphin = \bmun\,,
  \label{odeb}
\end{align}
where the matrices $A=(a_{ij})_{i,j=1}^n$ and $B=(b_{ij})_{i,j=1}^n$
are defined~by
\Beq
  a_{ij} := \iO \nabla\ej \cdot \nabla\ei
  \aand
  b_{ij} := \lambda_i \delta_{ij},
  \quad \hbox{for $i,j=1,\dots,n$},
  \non
\Eeq
and $\calB:\erren\to\erren$, $\calF:\erren\to\erren$ and $\calG:\erren\to\erren$ are the functions whose $i$th~components\last{, for $i=1,\dots,n$,} are defined~by
\begin{align}
  & \calB_i(s_1,\dots,s_n)
  := - \iO \Bigl\{ \beta''\Bigl( {\textstyle\somma j1n} s_j\ej \Bigr) \Bigl| \somma j1n s_j \nabla\ej \Bigr|^2 \, \ei \Bigr\}\,,
  \non
  \\
  & \calF_i(s_1,\dots,s_n)
  := \iO \Bigl\{ \beta'\Bigl( {\textstyle\somma j1n} s_j\ej \Bigr) \somma j1n s_j \lambdaj \ej \ei \Bigr\}\,,
  \non
  \\
  & \calG_i(s_1,\dots,s_n)
  := \iO \Bigl\{ (\gamma+g) \Bigl( {\textstyle\somma j1n} s_j\ej \Bigr) \ei \Bigr\} \,.
  \non
\end{align}	
Notice that these functions are well-defined, since the eigenfunctions are bounded \juerg{in view of the 
embedding $W \emb \Lx\infty$}.
For the same reason, one can differentiate under the integral sign,
so that they are (at~least) of class~$C^1$.
At this point, we eliminate $\bmun$ in \eqref{odea} by using \eqref{odeb} as definition of~$\bmun$\juerg{,
thus obtaining} a system of ordinary differential \revis{equations} in the only unknown~$\bphin$.
The local \Lip\ continuity just observed ensures that every Cauchy problem for this system 
has a unique maximal solution (which is even more regular than required).
This is the case if the initial condition for $\bphin$ is that derived from~\eqref{cauchyn},
so that the discrete problem \Pbln\ has a unique maximal solution \an{$\bphin$}
defined in \pier{some interval $[0,\Tn]$, with} $\Tn\in(0,T]$.
To complete the proof, we have to show that $\Tn=T$.
We do this by proving some a priori estimates \an{that are uniform with respect to $n$}.
Since the \juerg{established bounds prove to be} useful in the following,
we start to enumerate them.

\step
First a priori estimate

\an{To begin with, we recall the previous comment concerning the space $V_1$. We then test}
\eqref{priman} by~$1/|\Omega|\last{\in V_1}$ \an{to} obtain an ordinary differential equation for~$\phinbar$.
Namely, we have~that
\Beq
  \ddt \, \phinbar + \sigma\phinbar = \ubar
  \quad \juerg{\mbox{a.e. in $(0,T_n)$}} \,.
  \label{primanbar}
\Eeq
By also observing that $\phinbar(0)=\phizbar$, we immediately deduce~that
\Beq
  \norma\phinbar_{L^\infty(0,T_n)}
  \leq \Cstar \,,
  \label{primastiman}
\Eeq
where\juerg{, for further reference, we have used the special symbol $\Cstar$ instead of $c$}.

\step
Second a priori estimate

We recall the definition \eqref{defN} of~$\calN$ and~\eqref{calNvn},
and we \juerg{interpret} the functions $\phinbar$ and $\ubar$ as space-independent functions defined in~\juerg{$\Omega\times(0,T_n)$},
and thus \eqref{primanbar} as a partial differential equation.
\juerg{We multiply it by $v$, integrate over~$\Omega$, subtract the result from \eqref{priman}, and} test the \an{resulting equation} 
by $\calN(\phin-\phinbar)$.
At the same time, we test \eqref{secondan} by $\phin-\phinbar$.
Then\an{,} we add the resulting equalities to each other.
By recalling the properties \PropN\ of~$\calN$, which also produce a cancellation, \an{we} rearrang\an{e the terms}
\juerg{and infer that a.e. in $(0,T_n)$ it holds that}
\begin{align}
  & \ddt \, \normaVp{\phin-\phinbar}^2
  + \sigma \normaVp{\phin-\phinbar}^2
  \non
  \\
  & \quad {}
  + \iO |\Delta\phin|^2
  + \iO \beta'(\phin) |\nabla\phin|^2
  + \iO \gamma(\phin) (\phin-\phinbar)
  \non
  \\
  \separa
  & = \iO \last{ \beta'(\phin)\Delta\phin \, }(\phin-\phinbar)
  - (2\lambda-\nu) \iO \Delta\phin \, (\phin-\phinbar) 
  \non
  \\
  & \quad {}
  - \iO g(\phin) (\phin-\phinbar)
  + \iO u \, \calN(\phin-\phinbar)\,, 
  \label{persecondastima}
\end{align}
where in the last term we have written $u$ in place of $u-\ubar$ since $\calN(\phin-\phinbar)$ is orthogonal to \an{the space of constant functions $V_1$ to which $\ubar$ belongs.}
As for the \lhs, we notice that $\beta'$ is nonnegative, so that we just have to deal with the last term,
 \an{but this can be readily handled by using} that $\gamma$ is monotone, \an{ as remarked in \eqref{betagamma}. Namely, there exist}
some $\delta_0>0$ and $C_0>0$, \an{such that}
\Beq
  \gamma(s) (s-s_0)
  \geq \delta_0 |\gamma(s)| - C_0
  \quad \hbox{for every $s\in\erre$ and $s_0\in[-\Cstar,\Cstar]$\,,}
  \qquad
  \label{trickMZ}
\Eeq
where $\Cstar$ is the constant introduced in \eqref{primastiman}.
This inequality is similar to the one proved in \cite[Appendix, Prop. A.1]{MiZe}
(see also the detailed argument given in \cite[p.~908]{GiMiSchi} with a fixed $s_0$, 
which works in the present case with only minor \an{modifications}).
By applying it, we \juerg{find}~that
\Beq
  \iO \gamma(\phin) (\phin-\phinbar)
  \geq \delta_0 \iO |\gamma(\phin)|
  - c \quad\mbox{\juerg{a.e. in $(0,T_n)$}}\,.
  \juerg{\label{Fiete1}}
\Eeq
Now, we consider \an{the terms on the} \rhs\ \an{of \eqref{persecondastima}}.
We have~that
\begin{align}
  & \iO \last{ \beta'(\phin)\Delta\phin \, } (\phin-\phinbar)
  = - \iO \nabla\phin \cdot \bigl( \beta''(\phin) (\phin-\pier{\phinbar}) \nabla\phin + \beta'(\phin) \nabla\phin \bigr)
  \non
  \\
  & = - \iO \bigl( \beta''(\phin) - \beta''(\phinbar) \bigr) (\phin- \phinbar) |\nabla\phin|^2
  - \iO \bigl( \beta''(\phinbar) (\phin-\phinbar) + \beta'(\phin) \bigr) |\nabla\phin|^2
  \non
  \\
  & \leq - \iO \bigl( \beta''(\phinbar) (\phin-\phinbar) + \beta'(\phin) \bigr) |\nabla\phin|^2
  \,\quad\mbox{\juerg{a.e. in $(0,T_n)$}},
	\label{Fiete2}
\end{align}
where we have used the monotonicity of $\beta''$ implied by~\eqref{hpbeta}.
To estimate the last integral, \juerg{we perform almost everywhere in $\Omega\times (0,T_n)$ a 
Taylor expansion of $\beta'$ about $\phinbar$ \pier{using the integral remainder}, and then first account} for the sign of $\beta'''$ and then for~\eqref{primastiman}.
\pier{Then,} we have \juerg{almost everywhere in $\Omega\times (0,T_n)$} that
\begin{align}
  & - \bigl( \beta''(\phinbar) (\phin-\phinbar) + \beta'(\phin) \bigr)
  \non
  \\
  &
  = - \beta'(\phinbar) - 2 \beta''(\phinbar) (\phin-\phinbar) - 
 \pier{\Bigl[\int_0^1 \beta''' \big( \phinbar +s (\phin-\phinbar) \big) (1-s)\,ds\Bigr]}
   (\phin-\phinbar)^2
  \non
  \\
  & \leq - \beta'(\phinbar) - 2 \beta''(\phinbar) (\phin-\phinbar) 
  \leq \hat C \, (|\phin|+1) \,,
  \label{Fiete3}
\end{align}
where  we have \an{employed} the special symbol $\hat C$ instead of $c$, since now we use its precise value.
Now observe that by \eqref{hpbetaprimo} we can find some $s_1>0$ such that 
\begin{equation*}
\juerg{\hat C(|s|+1)\,\leq\,\frac 1 2 \beta'(s) \quad\mbox{ for }\,|s|>s_1.}
\end{equation*}
\juerg{At this point, we fix an arbitrary $t\in (0,T_n)$ for which \eqref{persecondastima} 
and \eqref{Fiete1}--\eqref{Fiete3},
written at the argument $t$,
hold true (we know that this is actually the case for almost every $t\in (0,T_n)$). 
We then put $\Omega_1:=\{x\in\Omega:\,|\phin(x,t)|>s_1\}$ and 
conclude that}
\begin{align}
  & \hat C \iO \!(|\phin\juerg{(t)}|+1) |\nabla\phin\juerg{(t)}|^2
  = \hat C \int_{\Omega_1} \!(|\phin\juerg{(t)}|+1) |\nabla\phin\juerg{(t)}|^2
  + \hat C \int_{\Omega\setminus\Omega_1} \!\!(|\phin\juerg{(t)}|+1) |\nabla\phin\juerg{(t)}|^2
  \non
  \\
  & \leq \frac 12 \iO \beta'(\phin\juerg{(t)}) |\nabla\phin\juerg{(t)}|^2
  + \hat C \, (s_1+1) \iO |\nabla\phin\juerg{(t)}|^2\,,
  \non
\end{align}
\an{where we also used that $\beta'$ \last{is nonnegative}.}
\juerg{For the last integral, we make use of the compactness inequality \eqref{compactbis}, which yields that}
\begin{align}
  & \hat C \, (s_1+1) \iO |\nabla\phin\juerg{(t)}|^2
  \leq \frac 14 \iO |\Delta\phin\juerg{(t)}|^2
  + c \, \normaVp{\phin\juerg{(t)}}^2
  \non
  \\
  & \leq \frac 14 \iO |\Delta\phin\juerg{(t)}|^2
  + c \, \normaVp{\phin\juerg{(t)}-\phinbar\juerg{(t)}}^2 + c \,,
  \non
\end{align}
\juerg{where $c$ is independent of both $t$ and $n$}. 
\an{Thus, we} deduce that
\begin{align}
  & - \iO \bigl( \beta''(\phinbar) (\phin-\phinbar) + \beta'(\phin) \bigr) |\nabla\phin|^2
  \leq c \iO (|\phin|+1) |\nabla\phin|^2
  \non
  \\
  & \leq \frac 12 \iO \beta'(\phin) |\nabla\phin|^2 
  + \frac 14 \iO |\Delta\phin|^2
  + c \, \normaVp{\phin-\phinbar}^2 + c \quad\mbox{\juerg{a.e. in \,$(0,T_n)$}}.
  \label{Fiete4}
\end{align}

With this, the estimation of the first term on the \rhs\ of \eqref{persecondastima} is completed.
Next, by  Young's inequality and the compactness inequality~\eqref{compactbis}, we have~that
\begin{align}
  & - (2\lambda-\nu) \iO \Delta\phin \, (\phin-\phinbar) 
  \leq \frac 18 \iO |\Delta\phin|^2
  + c \iO |\phin-\phinbar|^2
  \non
  \\
  & \leq \frac 14 \iO |\Delta\phin|^2
  + c \, \normaVp{\phin-\phinbar}^2 \quad\mbox{\juerg{a.e. in $(0,T_n)$}}.
  \label{Fiete5}
\end{align}

\juerg{To treat} the next term of~\eqref{persecondastima}\juerg{, we observe} that \eqref{gammag} ensures\an{, similar \juerg{to} above,} the existence of some $s_2>0$ such~that
\Beq
  |g(s)| \bigl( |s| + \Cstar \bigr) 
  \leq \frac {\delta_0}2 \, |\gamma(s)|
  \quad \hbox{whenever $|s|>s_2$}\,,
  \non
\Eeq
where $\Cstar$ is the constant appearing in \eqref{primastiman} and $\delta_0$ is the same as in~\eqref{trickMZ}.

\juerg{We now consider an arbitrary $t\in(0,\pier{T_n})$ for which \eqref{persecondastima} and \eqref{Fiete1}--\eqref{Fiete5}, evaluated at $t$, hold true
(this is the case for a.e. $t\in (0,\pier{T_n})$), and \last{put} $\Omega_2:=\{x\in\Omega: \,|\phin(x,t)|>s_2\}$. 
Then\pier{, it turns out that}
\begin{align}
  & - \iO g(\phin\juerg{(t)}) (\phin\juerg{(t)}-\phinbar\juerg{(t)}) \non \\ 
  &\leq \int_{\Omega_2} |g(\phin\juerg{(t)})| \, (|\phin\juerg{(t)}| + \Cstar) 
  + \int_{\Omega\setminus\Omega_2}  |g(\phin\juerg{(t)})| \, (|\phin\juerg{(t)}| + \Cstar) 
  \non
  \\
  & \leq \frac {\delta_0}2 \iO |\gamma(\phin\juerg{(t)})|
  + c \,,
  \non
\end{align}
\juerg{where, again, $c$ is independent of both $t$ and $n$.}
Finally, we easily see that
\Beq
  \iO u \, \calN(\phin-\phinbar)
  \leq c \, \norma{u }_{{*}}\, \normaV{\calN(\phin-\phinbar)}
  \leq c \, \normaVp{\phin-\phinbar}^2 + c \quad\juerg{\mbox{a.e. in }\,(0,T_n)},
  \non
\Eeq
where we have used the inequalities $\norma{v}_{{*}}\leq c\,\norma v_\infty$ for $v\in\Linfty$ \an{with $v=u$},
 and $\norma u_\infty\leq M$, \an{which we assumed in \eqref{hpdati}}.
At this point, we \an{go back to} \eqref{persecondastima}, \an{collect} the estimates just proved, and rearrange \last{the terms}.
Ignoring a nonnegative term on the \lhs, we deduce~that
\begin{align}
  & \ddt \, \normaVp{\phin-\phinbar}^2
  + \frac 12 \iO |\Delta\phin|^2
  + \frac 12 \iO \beta'(\phin) |\nabla\phin|^2
  + \frac {\delta_0}2 \iO \pier{|\gamma(\phin)|} 
  \non
  \\
  & \leq c \, \normaVp{\phin-\phinbar}^2 + c
  \quad \hbox{a.e.\ in $(0,\Tn)$} .
  \non
\end{align}
Now, we integrate over $(0,t)$ with an arbitrary $t\in(0,\Tn]$.
This integration \an{produces} the initial value~$\phin(0)$,
and \juerg{since} $\phin(0)$ is the $H$\an{-}\juerg{orthogonal} projection of $\phiz$ onto~$\Vn$ by~\eqref{cauchyn},
we have that 
\Beq
  \normaVp{\phin(0)}
  \leq c \, \norma{\phin(0)}
  \leq c \, \norma\phiz \,.
  \non
\Eeq
Thus, by applying the Gronwall lemma, we can \an{conclude that}
\Beq
  \norma\phin_{L^\infty(0,\Tn;\Vp)}
  \leq c \,.
  \label{secondastiman}
\Eeq
Now, since $\erren$ is finite dimensional, there holds \juerg{with some $c_n>0$} the inequality
\Beq
  |(s_1,\dots,s_n)| \leq c_n \, \Bigl\| {\textstyle\somma j1n s_j\ej} \Bigr\|_\Vp
  \quad \hbox{for every $(s_1,\dots,s_n)\in\erren$}.
  \non
\Eeq
Hence, \eqref{secondastiman} implies that
\Beq
  \norma\bphin_{L^\infty(0,\Tn;\erren)} \leq c_n\,, 
  \non
\Eeq
that is, that $\bphin$ is \last{uniformly} bounded.
Since $\bphin$ is also maximal, it has to be global.
Thus, $\Tn=T$, and the proof is complete.}
\Edim

The rest of the present subsection is devoted to the existence and stability part of Theorem~\ref{Wellposedness}.
To this end, we first upgrade and improve the estimates already obtained, and then perform further ones.
The first estimate \eqref{primastiman} yields an $L^\infty$ bound for~$\phinbar$.
However, \juerg{we immediately find} from \eqref{primanbar} that the derivative of $\phinbar$ is bounded as well.
Hence, we have~that
\Beq
  \norma\phinbar_{W^{1,\infty}(0,T)}
  \leq c \,.
  \label{primastima}
\Eeq
Next, the above argument involving Gronwall's lemma, \juerg{which now can be applied} for arbitrary $t\in(0,T]$, 
leads \an{to} \juerg{the estimate}
\Beq
  \norma{\phin-\phinbar}_{\L\infty\Vp}
  + \norma{\Delta\phin}_{\L2H}
  + \norma{\beta'(\phin)|\nabla\phin|^2}_{\LQ1}
  + \norma{\gamma(\phin)}_{\LQ1}
  \leq c \,.
  \label{secondastima}
\Eeq
Thus, on account of \eqref{primastima}, \pier{we realize that $\phin$ is} bounded in $\L\infty\Vp$\last{, and by} also \pier{invoking the} elliptic regularity theory, we conclude~that
\Beq
  \norma\phin_{\L\infty\Vp\an{\cap \L2W}} \leq c \,.
  \label{dasecondastima}
\Eeq

\step
Third a priori estimate

\an{We now} test \eqref{secondan} \an{by~$1 \in V_1$}, \juerg{obtaining}~that
\Beq
  |\Omega| \, \munbar
  = - \iO \beta'(\phin) \Delta\phin
  + \iO \bigl( \gamma(\phin) + g(\phin) \bigr) \juerg{\quad\mbox{a.e. in $(0,T)$}}.
  \label{munbar}
\Eeq
Now, by also accounting for \an{the growth condition in} \eqref{hpbetasecondo}, we \an{infer} that
\begin{align}
  & \Bigl| - \iO \beta'(\phin) \Delta\phin \Bigr|
  = \Bigl| \iO \nabla\phin \cdot \nabla\beta'(\phin) \Bigr|
  \non
  \\
  & \leq \iO |\beta''(\phin)| \, |\nabla\phin|^2
  \leq c \iO (\beta'(\phin)+1) |\nabla\phin|^2  \juerg{\quad\mbox{a.e. in $(0,T)$}}.
  \label{auxterza}
\end{align}
On the other hand, the definition \eqref{defg} of~$g$, \eqref{gammag}, and the continuity of the \an{involved functions}, imply~that
\Beq
  |g(s)| \leq c \, (|\gamma(s)|+1)
  \quad \hbox{for every $s\in\erre$}.
  \label{stimag}
\Eeq
Hence, recalling \accorpa{secondastima}{dasecondastima}\an{,} we conclude that
\Beq
  \norma\munbar_{L^1(0,T)} \leq c \,.
  \label{terzastima}
\Eeq

\step
Fourth a priori estimate

We test \eqref{priman} by $\mun$ and have \aet~that
\Beq
  \iO \dt\phin \, \mun
  + \iO |\nabla\mun|^2
  = \iO (u-\sigma\phin) \mun\,.
  \label{Fiete6}
\Eeq
Next, we rewrite \eqref{secondan} in terms of the original potential, \last{as}
\begin{align}
  & \iO \Delta\phin \, \Delta v
  - \iO \Delta f(\phin) \, v
  - \iO f'(\phin) \Delta\phin \, v
  \non
  \\
  & \quad {}
  - \nu \iO \Delta\phin \, v
  + \iO f(\phin) \, f'(\phin) \, v
  + \nu \iO f(\phin) \, v
  \non
  \\
  & 
  = \iO \mun \, v
  \quad \hbox{\aet, for every $v\in\Vn$}\,,
  \non
\end{align}
and test it by~$\dt\phin$.
Moreover, we perform an integration by parts on the \lhs,
\juerg{where we recall} that both $\phin$ and $\dt\phin$ are $W$-valued.
We obtain~that
\begin{align}
  & \iO (-\Delta\phin) \, (-\Delta\dt\phin)
  - \iO f(\phin) \Delta\dt\phin
  - \iO f'(\phin) \Delta\phin \, \dt\phin
  \non
  \\
  & \quad {}
 - \nu \iO \Delta\phin \, \dt\phin + \iO f(\phin) \, f'(\phin) \, \dt\phin
  + \nu \iO f(\phin) \, \dt\phin
  \non
  \\
  & = \iO \mun \, \dt\phin \juerg{\quad\mbox{a.e. in $(0,T)$}},
  \label{Fiete7}
\end{align}
and recognize the time derivative of $\calE(\phin)$ (see \eqref{defE}) on the \lhs.
Hence, by adding \juerg{the equalities \eqref{Fiete6} and \eqref{Fiete7}} to each other\an{, rearranging,} and noting 
an obvious cancellation,
we deduce~that
\Beq
  \an{ 
  \ddt \, \calE(\phin)
  + \iO |\nabla\mun|^2
  }
  = \iO (u-\sigma\phin) \mun \juerg{\quad\mbox{a.e. in $(0,T)$}}.
  \non
\Eeq
Integration over $(0,t)$\last{, for an arbitrary $t \in (0,T]$,} then yields~that\an{, recalling} the definition \eqref{defQt} of~$Q_t$\an{,}
\begin{align}
  &
  \frac 12 \iO |{-\Delta\phin(t)} + f(\phin(t))|^2
  + \frac \nu 2 \iO |\nabla\phin(t)|^2 
  + \nu \iO F(\phin(t)) 
  \an{
  +  \intQt |\nabla\mun|^2
  }
  \non
  \\
  & = \frac 12 \iO |{-\Delta\phin(0)} + f(\phin(0))|^2
  + \frac \nu 2 \iO |\nabla\phin(0)|^2 
  + \nu \iO F(\phin(0)) 
  \non
  \\
  & \quad {}
  + \intQt (u-\sigma\phin)\mun 
  \quad \hbox{for every $t\in[0,T]$}.
  \label{perquartastima}
\end{align}
The \lhs\ is bounded from below (see \eqref{bddbelow}).
As for the terms involving the initial value, 
we recall that $\phiz\in W$ and that $\normaW{\phin(0)}$ is bounded
(see \eqref{hpdati} and Remark~\ref{RemVn}).
Thus, also $\norma{\phin(0)}_\infty$ is bounded due to the continuous embedding $W\emb\Linfty$,
\an{so that} all of the terms involving the initial value are bounded.
Let us estimate the last term \an{on the \rhs\ of}~\eqref{perquartastima}.
\an {It holds that}
\Beq
  \intQt (u-\sigma\phin)\mun
  = \intQt (u-\sigma\phin) (\mun-\munbar)
  + \intQt (u-\sigma\phin)\munbar\,,
  \non
\Eeq
where the Young and Poincar\'e inequalities imply that
\Beq
  \intQt (u-\sigma\phin) (\mun-\munbar)
  \leq \frac 12 \intQt |\nabla\mun|^2
  + c \intQt ( |u|^2 + |\phin|^2 \bigr)
  \leq \frac 12 \intQt |\nabla\mun|^2 
  + c \,.
  \non
\Eeq
In the last inequality, we \juerg{made use of} \eqref{hpdati} and \eqref{terzastima},
which also justify the \an{computations to follow}.
Indeed, owing to~\eqref{dasecondastima} \pier{as well}, we can estimate the remaining term this~way:
\begin{align}
  & \intQt (u-\sigma\phin)\munbar
  \leq c \ioT \norma{u\an{(s)}-\sigma\phin(s)}_* \, \normaV{\munbar(s)} \, ds
  \non
  \\
  & {} = c \ioT \norma{u\an{(s)}-\sigma\phin(s)}_* \, |\munbar(s)| \, ds
  \leq c \, \norma{u-\an{\sigma\phin}}_{\L\infty\Vp} \, \norma\munbar_{L^1(0,T)}
  \leq c \,.
  \non
\end{align}
Collecting \an{the above \juerg{computations}}, we conclude that
\begin{align}
  & \norma{-\Delta\phin+f(\phin)}_{\L\infty H}
  + \norma{\nabla\phin}_{\L\infty H}
  \non
  \\
  & \quad {}
  + \norma{F(\phin)}_{\L\infty\Luno}
  \an{
  + \norma{\nabla\mun}_{\L2H}
  }
  \leq c \,.
  \non
\end{align}
By a standard argument regarding the elliptic operator $-\Delta+f(\cdot)$\an{,}
and accounting for \eqref{secondastima} once more, we deduce~that
\begin{align}
  & \norma\phin_{\L\infty W}
  + \norma{f(\phin)}_{\L\infty H}
  + \norma{F(\phin)}_{\L\infty\Luno}
  \an{
  + \norma{\nabla\mun}_{\L2H}
  }
  \leq c \,.
  \label{quartastima}
\end{align}

\step
Fifth a priori estimate

Thanks to the continuous embedding $W\emb\Linfty$, we have that $\L\infty W\emb\LQ\infty$.
Therefore, \juerg{also using the regularity of the nonlinearities,
we deduce from}  \eqref{quartastima} that
\Beq
  \norma\phin_\infty \leq c \,, \quad
  \norma{\beta'(\phin)}_\infty \leq c\,,
  \aand
  \norma{\gamma(\phin)}_\infty \leq c\,, \juerg{\quad\mbox{a.e. in $(0,T)$}}\,.
  \non
\Eeq
Hence, coming back to \eqref{munbar} and \eqref{stimag}\an{,} we conclude that
\Beq
  \norma\munbar_{L^\infty(0,T)} \leq c\,,
  \non
\Eeq
\an{\juerg{which, together} with \eqref{quartastima} and \Poincare's inequality, produces}
\Beq
  \norma\mun_{\L2V} \leq c \,.
  \label{quintastima}
\Eeq

\step 
Sixth a priori estimate

We take any $v\in\L2V$ and define $\vn\in\L2\Vn$ by setting
$\vn(t):=\Pn v(t)$ \aat\an{.}
Then, we test \eqref{priman} by~$\vn$.
Noting that $\dt\phin(t)\in\Vn$ \aat, and owing to the above estimates \last{along with} Remark~\ref{RemVn}, we \last{obtain}~that
\begin{align}
  & \intQ \dt\phin \, v
  = \intQ \dt\phin \, \vn
  = \intQ (u-\sigma\phin) \vn 
  - \an{\intQ} \nabla\mun \cdot \nabla\vn
  \non
  \\
  & \leq c \, \norma\vn_{\L2V}
  \leq c \, \norma v_{\L2V}\,,
  \non
\end{align}
\an{from which it \last{is readily seen} that}
\Beq
  \norma{\dt\phin}_{\L2\Vp} \leq c \,.
  \label{sestastima}
\Eeq

\step
Conclusion of the existence proof

Eventually, by letting $n$ tend to infinity, we find a solution $(\phi,\pier{\mu})$ 
to the equations \eqref{prima} and \eqref{secondabis} \juerg{that satifies}  the initial condition~\eqref{cauchy}.
\juerg{Recalling all the above estimates, and invoking} well-known compactness results
(for the strong compactness see, e.g., \cite[Sect.~8, Cor.~4]{Simon}),
we have that\an{, as $n\to\infty$,}
\begin{align}
  & \phin \to \phi
  \quad \hbox{weakly star in $\H1\Vp\cap\L\infty W$,}
  \non
  \\ 
  & \quad \hbox{strongly in $\C0V$ and uniformly in $Q$}\,,
  \label{convphin}
  \\
  & \mun \to \mu
  \quad \hbox{weakly in $\L2V$}\,,
  \label{convmun}
\end{align}
for some pair $(\phi,\mu)$, and at least for a (nonrelabeled) subsequence.
Notice that also the inequality obtained from \eqref{stability} by ignoring the part regarding $w$ is satisfied, 
due to the \an{semicontinuity} of norms.

By \eqref{convphin}, $\{\phin(0)\}$ converges to $\phi(0)$ strongly in~$V$.
On the other hand, since $\phin(0)$ is the $H$-projection of~$\phiz$ onto~$\Vn$, it converges to $\phiz$ \an{strongly} in~$H$.
Hence, the initial condition \eqref{cauchy} is satisfied.
We show that the equations \eqref{prima} and \eqref{secondabis} are satisfied as well.
Since $\{\phin\}$ is uniformly bounded in $\LQ\infty$, and \juerg{since} the \an{involved nonlinearities} in the original system are smooth
\juerg{and}
thus \Lip\ continuous on bounded sets, we can deduce~that
\Beq
  \beta^{(i)}(\phin) \to \beta^{(i)}(\phi)
  \quad \hbox{uniformly in $Q$,\quad for $i=0,1,2$}.
  \label{convbeta}
\Eeq
At this point, we write the time-integrated version of \eqref{priman} and \eqref{secondan}, namely
\begin{align}
  & \intQ \dt\phin \, \vn 
  \an{
  + \intQ \nabla\mun \cdot \nabla\vn
  + \sigma \intQ \phin \vn
  }
  = \intQ u \vn
  \quad \hbox{for every $\vn\in\L2\Vn$}\,,
  \label{intpriman}
  \\
  \separa
  & \intQ \Delta\phin \, \Delta\vn
  - \intQ \Delta\beta(\phin) \, \vn
  - \intQ \beta'(\phin) \Delta\phin \, \vn
  \non
  \\
  & \quad {}
  + \intQ \gamma(\phin) \, \vn
  + (2\lambda-\nu) \an{\intQ} \Delta\phin \, \vn
  + \intQ g(\phin) \, \vn
  \non
  \\
  & = \intQ \mun \, \vn
  \quad \hbox{for every $\vn\in\L2\Vn$}\,.
  \label{intsecondan}
\end{align}
Our aim is \juerg{to take the limit as $n\to\infty$} in these equation\an{s} with $V$\an{-} or $W$-valued test functions.
We fix $v\in\L2V$ and take as $\vn$ in \eqref{intpriman} the function defined by setting $\vn(t):=\Pn v(t)$ \aat.
Then, $\{\vn\}$ converges to $v$ strongly in $\L2V$ by Remark~\ref{RemVn}, and we have~that
\Beq
  \intQ \dt\pier\phi \, v
  \an{{}
  + \intQ \nabla\pier\mu \cdot \nabla v
  + \sigma \intQ \pier\phi\, v
  }
  = \intQ u v \,.
  \non
\Eeq
Since $v$ is arbitrary in $\L2V$, this is equivalent to~\eqref{prima}.
Similarly, given an arbitrary $v\in\L2W$, we test \eqref{intsecondan} by the function $\vn$ defined as before.
Now, $\{\vn\}$~converges to $v$ strongly in~$\L2W$.
To take the limit in \eqref{intsecondan}, it is convenient to transform one of the terms \juerg{as follows:}
\begin{align}
  - \intQ \Delta\beta(\phin) \, \vn
  = \intQ \beta'(\phin) \nabla\phin \cdot \nabla\vn \,.
  \non
\end{align}
Then we obtain~that
\begin{align}
  & \intQ \Delta\phi \, \Delta v
  + \intQ \beta'(\phi) \nabla\phi \cdot \nabla v
  - \intQ \beta'(\phi) \Delta\phi \, v
  \non
  \\
  & \quad {}
  + \intQ \gamma(\phi) \, v
  + (2\lambda-\nu) \iO \Delta\phi \, v
  + \intQ g(\phi) \, v
  = \intQ \mu \, v \,.
  \label{limsecondan}
\end{align}
Now, we have that $\{\phin\}$ converges to $\phi$ strongly in $\C0{\Wx{1,4}}$
thanks to the \an{just quoted} strong compactness result \an{and} the compactness of the embedding $W\emb\Wx{1,4}$.
It follows that $\{\nabla\phin\}$ converges to $\nabla\phi$ strongly in $\C0{\Lx4}$,
and \juerg{$\{|\nabla\phin|^2$\}} converges to $|\nabla\phi|^2$ strongly in $\C0H$.
Since
\Beq
  \Delta\beta(\phin) 
  = \beta''(\phin) |\nabla\phin|^2 + \beta'(\phin)\Delta\phin\,,
  \non
\Eeq
we deduce that $\{\Delta\beta(\phin)\}$ strongly converges (at least) in $\L2H$\last{, as $n\to \infty$,}
\an{and that its limit} must be $\Delta\beta(\phi)$, due to~\eqref{convbeta}.
Hence we can write
\Beq
  \intQ \beta'(\phi) \nabla\phi \cdot \nabla v
  = - \intQ \Delta\beta(\phi) \, v 
  \non
\Eeq
in \eqref{limsecondan}.
Since $v$ is arbitrary in $\L2W$,
we \juerg{infer} that \eqref{limsecondan} is equivalent to~\eqref{secondabis}, \juerg{thus}
\an{concluding the proof}.


\subsection{Uniqueness and continuous dependence}
\label{CONTDEP}

This subsection is devoted to the proof of the uniqueness of the solution to the problem~\Pbl\ 
and the continuous dependence estimate \eqref{contdep}.
Thanks to Proposition~\ref{Pblbis}, we can deal with the alternative \an{formulation in}~\eqref{pblbis}.
\an{The outline of our} strategy is \juerg{as follows}:\,
$(i)$~\an{prove} an inequality in the direction of~\eqref{contdep} for arbitrary pairs of solutions;\,
$(ii)$~\an{derive} uniqueness;\,
$(iii)$~\an{complete} the proof of~\eqref{contdep}. 
\an{First, let us} fix two functions $u_i\in\LQ\infty$ \juerg{with $\norma{u_i}_\infty\leq M$, $i=1,2$}, 
consider arbitrary corresponding solutions $(\phi_i,\mu_i)$ \an{as given by Theorem \ref{Wellposedness}, and} set for convenience
\Beq
  u := u_1 - u_2 \,, \quad
  \phi := \phi_1 - \phi_2\,,
  \aand
  \mu := \mu_1 - \mu_2 \,.
  \non
\Eeq
We first write \eqref{prima} for both solutions and take the difference.
We thus \an{find}~that
\Beq
  \< \dt\phi , v >_\an{V}
  \an{{}
  + \iO \nabla\mu \cdot \nabla v
  + \sigma \iO \phi v
  }
  = \iO u v
  \quad \hbox{\aet, for every $v\in V$}\,. 
  \label{diffprima}
\Eeq
We \an{test} this equation by $1/|\Omega|$,
\an{obtaining an identity} for $\phibar$\juerg{, which we then multiply by any $v\in V$ to obtain the equation}
\Beq
  \< \dt\phibar , v >_\an{V}
  + \sigma \iO \phibar v
  = \iO \ubar v
  \quad \hbox{\aet, for every $v\in V$} .
  \label{diffprimabar}
\Eeq
Next, we take the difference between \eqref{diffprima} and \eqref{diffprimabar}
and test it by $\calN(\phi-\phibar)$.
Recalling the properties \PropN, we have \aet~that
\Beq
  \frac 12 \, \ddt \, \normaVp{\phi-\phibar}^2
  \an{{}
  + \iO \mu (\phi-\phibar)
  + \sigma \, \normaVp{\phi-\phibar}^2
  }
  = \iO (u-\ubar) \calN(\phi-\phibar)  \,.
  \label{testdiffprima}
\Eeq
At the same time, we write \eqref{secondabis} for both solutions and test the difference by $\phi-\phibar$.
\an{Upon} rearranging, we obtain~that
\begin{align}
  & \iO |\Delta\phi|^2
   = - \iO \nabla \bigl(
      \beta(\phi_1) - \beta(\phi_2)
    \bigr) \cdot \nabla\phi
  \non
  \\
  & \quad {}
  + \iO \bigl(
    \beta'(\phi_1) \Delta\phi_1
    - \beta'(\phi_2) \Delta\phi_2
  \bigr) (\phi-\phibar)
  \pier{{}- \iO \bigl(
    \gamma(\phi_1)
    - \gamma(\phi_2) 
  \bigr) (\phi - \phibar)}
  \non
  \\  
  & \quad {}
  + (2\lambda-\nu) \iO |\nabla\phi|^2 
  \pier{{}- \iO \bigl(
    g(\phi_1)
    - g(\phi_2)
  \bigr) (\phi - \phibar)}
  + \iO \mu (\phi-\phibar) \,.
  \label{testdiffsecondabis}
\end{align}
Now, we estimate the \last{integrals on the} \rhs.
For a while, we allow the values of $c$ to depend on the solutions we have fixed.
In particular, \juerg{since} $\phi_1$ and~$\phi_2$ are bounded,
we \juerg{can make use of} the \Lip\ continuity of the nonlinearities on their range,
the \Lip\ constant being dependent on these solutions.
\juerg{To begin with, we have for the \rhs\ of \eqref{testdiffprima} that}
\begin{align}
& \juerg{\iO (u-\ubar) \calN(\phi-\phibar)\,=\,\iO u \calN(\phi-\phibar) 
\,\le\, \|\calN (\phi-\phibar)\|_V\, \|u\|_{V^*}} \non \\
&\juerg{\le\, c\,\|\phi-\phibar\|_*\, \|u\|_* \,\le\, \frac\sigma 2 \, \normaVp{\phi-\phibar}^2
  + c \, \normaVp u^2\,.}
\non	
\end{align}

Next, we deal with the \rhs\ of \eqref{testdiffsecondabis}.
To estimate the first term, we recall that $\beta'$ is nonnegative, 
the continuity of the embedding $W\emb\Wx{1,4}$, and the regularity \eqref{regphi} for~$\phi_1$.
Moreover, we owe to the \Lip\ continuity of $\beta'$ and the \Holder\ inequality, 
and apply the compactness inequality~\eqref{compactbis}.
For every $\delta>0$, \an{it holds}~that
\begin{align}
  & - \iO \nabla \bigl(
    \beta(\phi_1) - \beta(\phi_2)
  \bigr) \cdot \nabla\phi
  = \an{-}\iO \bigl(
    \beta'(\phi_1)\nabla\phi_1
    - \beta'(\phi_2)\nabla\phi_2
  \bigr) \cdot \nabla\phi
  \non
  \\
  & = \an{-}\iO \bigl(
    \beta'(\phi_1) - \beta'(\phi_2)
  \bigr) \nabla\phi_1 \cdot \nabla\phi
  - \iO \beta'(\phi_2) |\nabla\phi|^2
  \leq c \iO |\phi| \, \an{|}\nabla\phi_1| \, |\nabla\phi|
  \non
  \\
  & \leq c \, \norma\phi_4 \, \norma{\nabla\phi_1}_4 \, \norma{\nabla\phi}
  \leq c \, \normaV\phi^2
  \leq \delta \iO |\Delta\phi|^2
  + \cdelta \, \normaVp\phi^2 \,.
  \non
\end{align}
Now observe that the procedure used to derive \eqref{primanbar}
can be applied to $\phi$ and~$u$, 
so that the pointwise values of $\phibar$ can be estimated by the $L^2$ norm of~$\ubar$,
and thus also by the norm of $u$ in~$\L2\Vp$. \juerg{We therefore have}
$$
\juerg{|\phibar|\,\le\,c\,\|u\|_{\L2\Vp} \quad\mbox{a.e. in }\,(0,T),}
$$
\juerg{whence we conclude that}
\Beq
  \normaVp\phi^2
  \leq 2 \, \normaVp{\phi-\phibar}^2 + 2 \, \normaVp\phibar^2
  \leq 2 \, \normaVp{\phi-\phibar}^2 + c \, |\phibar|^2
  \leq 2 \, \normaVp{\phi-\phibar}^2 + c \, \norma u_{\L2\Vp}^2 \,.
  \non
\Eeq
Therefore,  \last{we find that}
\Beq
  - \iO \nabla \bigl(
    \beta(\phi_1) - \beta(\phi_2)
  \bigr) \cdot \nabla\phi
  \leq \delta \iO |\Delta\phi|^2
  + \cdelta \, \normaVp{\phi-\phibar}^2
  + \cdelta \, \norma u_{\L2\Vp}^2 \,.
  \non
\Eeq
\an{Mimicking} \juerg{this} \an{argument}, we \juerg{find for} the next term that
\begin{align}
  & \iO \bigl(
    \beta'(\phi_1) \Delta\phi_1
    - \beta'(\phi_2) \Delta\phi_2
  \bigr) (\phi-\phibar)
  \non
  \\
  & = \iO \bigl( \beta'(\phi_1) - \beta'(\phi_2) \bigr) \Delta\phi_1 \, (\phi-\phibar)
  + \iO \beta'(\phi_2) \Delta\phi \, (\phi-\phibar)
  \non
  \\
  & \leq c \, \norma\phi_4 \, \norma{\Delta\phi_1} \, \norma{\phi-\phibar}_4
  + c \, \norma{\Delta\phi} \, \norma{\phi-\phibar}
  \non
  \\
  & \leq c \, \normaV\phi \, \normaV{\phi-\phibar}
  + c \, \norma{\Delta\phi} \, \norma{\phi-\phibar}
  \non
  \\
  & \leq c \, \normaV{\phi-\phibar}^2
  + c \, \norma\ubar_{\pier{L^2(0,T)}}^2
  + \delta \iO |\Delta\phi|^2
  + \cdelta \, \norma{\phi-\phibar}^2
  \non
  \\
  & \leq 3\delta \iO |\Delta\phi|^2
  + \cdelta \, \normaVp{\phi-\phibar}^2
  + c \, \norma u_{\L2\Vp}^2 \,.
  \non
\end{align}
Since $\gamma$ and $g$ are \Lip\ continuous on the range of the solutions \juerg{under consideration},
the sum $S$ of the next \pier{three} terms can be treated as before using the compactness inequality\juerg{, namely}
\Beq
  S \pier{{}\leq \norma{\phi}^2\norma{\phi-\phibar}^2 + \normaV{\phi-\phibar}^2}
  \leq \pier{2} \delta \iO |\Delta\phi|^2 + \cdelta \, \normaVp{\phi-\phibar}^2 \pier{{}+ c \, \norma u_{\L2\Vp}^2}\,.
  \non
\Eeq
At this point, we add \eqref{testdiffprima} and \eqref{testdiffsecondabis} to each other
and notice \an{the} cancellation \an{of the terms involving the chemical potentials}.
Then, we account for all of the above estimates, choose $\delta$ small enough,
integrate with respect to time, and apply the Gronwall lemma \an{to} conclude~that
\Beq
  \norma\phi_{\L\infty\Vp}
  + \norma{\Delta\phi}_{\L2H}
  \leq c \, \norma u_{\L2\Vp} \,.
  \label{provvisoria}
\Eeq
In this estimate, the constant $c$ depends on the solutions we have considered.
Nevertheless, by applying it in the case $u_1=u_2$, we derive that $\phi_1=\phi_2$.
Then, by also considering the third components $w_i$ of the solutions to the original problem,
and recalling first \eqref{terza} and then~\eqref{seconda},
we deduce that $w_1=w_2$ and $\mu_1=\mu_2$.
This proves the uniqueness part of Theorem~\ref{Wellposedness}.

Now, we come back to \eqref{provvisoria} and observe that the uniqueness just established
implies that the solutions $\phi_i$\an{, for $i=1,2$,} at hand must coincide with those constructed in our existence proof.
Therefore, their norms we have considered in this proof are \an{uniformly} bounded as specified in the stability inequality.
We conclude that the constant $c$ that appears in \eqref{provvisoria} 
can be estimated in terms of $\Omega$, $T$, the structure of the system,
the initial datum $\phiz$ and~$M$.

So, from now on, the symbol $c$ denotes constants that just depend 
on $\Omega$, $T$, the structure of the system, the initial datum~$\phiz$, and~$M$, as before,
since the norms of the solutions involved in our argument can be controlled by the stability estimate~\eqref{stability}.
In particular, this \juerg{holds} for the estimate
\Beq
  \norma\phi_{\L\infty\Vp \an{\cap}\L2W}
  \leq c \, \norma u_{\L2\Vp}\,, 
  \label{primacontdep}
\Eeq
which is obtained from \eqref{provvisoria} on account of elliptic regularity.

\an{We are now going to continue the stability analysis of \pier{\eqref{prima}--\eqref{cauchy}} 
with respect to the control variable. This will be a crucial result
\juerg{in the derivation of} first-order optimality conditions for the associated optimal control problem.
Thus,} we write \eqref{prima} for both solutions and test the difference by~$\phi$.
At the same time, recalling \eqref{regDeltaphi}, we write \eqref{secondabis} for both solutions
and test the difference by $-\Delta\phi$.
Then we add the \juerg{resulting equalities} to each other and rearrange,
\an{noticing that a cancellation occurs}. \pier{By this procedure, we} \an{obtain} \aet~that
\begin{align}
  & \frac 12 \, \ddt \iO |\phi|^2
  + \sigma \iO |\phi|^2
  + \iO |\nabla\Delta\phi|^2
  \non
  \\
  & = \iO u \phi
  - \iO \bigl( \Delta\beta(\phi_1) - \Delta\beta(\phi_2) \bigr) \, \Delta\phi
  \non
  \\
  & \quad {}
  - \iO \bigl( \beta'(\phi_1)\Delta\phi_1 - \beta'(\phi_2)\Delta\phi_2 \bigr) \Delta\phi
  + \iO \bigl( \gamma(\phi_1) - \gamma(\phi_2) \bigr) \, \Delta\phi
  \non
  \\
  & \quad {}
  + (2\lambda-\nu) \iO |\Delta\phi|^2
  + \iO \bigl( g(\phi_1) - g(\phi_2) \bigr) \, \Delta\phi \,,
  \label{persecondacontdep}
\end{align}
and we have to estimate the \an{terms on the}~\rhs.
The first \an{one} is trivial, and the second integral is estimated as follows:
\begin{align}
  & - \iO \bigl( \Delta\beta(\phi_1) - \Delta\beta(\phi_2) \bigr) \, \Delta\phi
  = \iO \bigl( \beta'(\phi_1) \nabla\phi_1 - \beta'(\phi_2) \nabla\phi_2 \bigr) \cdot \nabla\Delta\phi
  \non
  \\
  & \leq \delta \iO |\nabla\Delta\phi|^2
  + \cdelta \iO |\beta'(\phi_1) \nabla\phi_1 - \beta'(\phi_2) \nabla\phi_2|^2 \,\last{,}
  \non
\end{align}
\last{where $\delta$ is an arbitrary positive constant to be chosen later on.}
Now, by accounting for the \Lip\ continuity of~$\beta'$, \Holder's inequality,
the regularity of~$\phi_1$, and the continuity of the embedding $V\emb\Lx4$,
we can~\an{bound the last integral in the following form:}
\begin{align}
  & \iO |\beta'(\phi_1) \nabla\phi_1 - \beta'(\phi_2) \nabla\phi_2|^2
  \leq 2 \iO |\beta'(\phi_1) - \beta'(\phi_2)|^2 \, |\nabla\phi_1|^2
  + 2 \iO |\beta'(\phi_2)|^2 \, |\nabla\phi|^2
  \non
  \\
  & \leq c \, \norma\phi_4^2 \, \norma{\nabla\phi_1}_4^2
  + c \, \normaV\phi^2
  \leq c \, \normaV\phi^2 \,.
  \non
\end{align}
\an{For the third term \juerg{on the \rhs\ of \eqref{persecondacontdep}}, we proceed}
with a similar argument, \last{obtaining that}
\begin{align}
  & - \iO \bigl( \beta'(\phi_1)\Delta\phi_1 - \beta'(\phi_2)\Delta\phi_2 \bigr) \Delta\phi
  = - \iO \an{\big(} \bigl( \beta'(\phi_1) - \beta'(\phi_2) \bigr) \Delta\phi_1 + \beta'(\phi_2) \Delta\phi \an{\big)} \Delta\phi
  \non
  \\
  & \leq c \, \norma\phi_4 \, \norma{\Delta\phi_1} \, \norma{\Delta\phi}_4
  + c \, \iO |\Delta\phi|^2
  \leq c \, \normaV\phi^2 
  + \delta \iO |\nabla\Delta\phi|^2
  + \cdelta \, \normaW\phi^2 \,.
  \non
\end{align}
Finally, the last three terms of \eqref{persecondacontdep} can be treated in a \pier{straightforward} way.
At this point, we integrate \eqref{persecondacontdep} with respect to time,
account for all of the estimates we have established, choose $\delta$ small enough,
and apply \eqref{primacontdep} to conclude that
\Beq
  \norma\phi_{\C0H\cap\L2{\Hx3}} \leq c \, \norma u_{\L2\Vp} \,.
  \label{secondacontdep}
\Eeq

Now, we recall Remark~\ref{Pblter}\an{,} write \eqref{prima} for both solutions,
and test the difference by $-\Delta\phi$.
At the same time, we write \pier{\eqref{secondabisPde}} for both solutions, 
multiply the difference by $\Delta^2\phi$, and integrate over~$\Omega$.
Then, we add the \juerg{resulting equalities} to each other.
Since a cancellation occurs in the terms involving~$\mu$,
we have \aet~that
\begin{align}
  & \frac 12 \, \ddt \iO |\nabla\phi|^2
  + \sigma \iO |\nabla\phi|^2
  + \iO |\Delta^2\phi|^2
  \non
  \\
  & = \iO u (-\Delta\phi)
  + \iO \bigl( \Delta\beta(\phi_1) - \Delta\beta(\phi_2) \bigr) \, \Delta^2\phi
  \non
  \\
  & \quad {}
  + \iO \bigl( \beta'(\phi_1)\Delta\phi_1 - \beta'(\phi_2)\Delta\phi_2 \bigr) \, \Delta^2\phi
  - \iO \bigl( \gamma(\phi_1) - \gamma(\phi_2) \bigr) \, \Delta^2\phi
  \non
  \\
  & \quad {}
  - (2\lambda-\nu) \iO \Delta\phi \, \Delta^2\phi
  - \iO \bigl( g(\phi_1) - g(\phi_2) \bigr) \Delta^2\phi \,.
  \label{perterzacontdep}
\end{align}
As for the \rhs, we immediately have that
\Beq
  \iO u (-\Delta\phi)
  \leq c \, \normaVp u \, \normaV{\Delta\phi}
  \leq \normaVp u^2 + c \, \norma\phi_{\Hx3}^2\,,
  \non
\Eeq
while the next \pier{terms need} some treatment.
Since $\Delta\beta(\phi_i)=\beta'(\phi_i)\Delta\phi_i+\beta''(\phi_i)|\nabla\phi_i|^2$, \an{for $i=1,2$,}
we \pier{can handle the second and third term on the \rhs\ as}
\begin{align}
  & \iO \bigl( \Delta\beta(\phi_1) - \Delta\beta(\phi_2) \bigr) \, \Delta^2\phi
   \pier{{}+ \iO \bigl( \beta'(\phi_1)\Delta\phi_1 - \beta'(\phi_2)\Delta\phi_2 \bigr) \, \Delta^2\phi}
  \non
  \\
  &  
  \last{{}={}}2 \iO \bigl( \beta'(\phi_1)\Delta\phi_1 - \beta'(\phi_2)\Delta\phi_2 \bigr) \, \Delta^2\phi
  + \iO \bigl( \beta''(\phi_1)|\nabla\phi_1|^2 - \beta''(\phi_2)|\nabla\phi_2|^2 \bigr) \, \Delta^2\phi
  \non
  \\
  & \leq \delta \iO |\Delta^2\phi|^2
  + \cdelta \iO |\beta'(\phi_1)\Delta\phi_1 - \beta'(\phi_2)\Delta\phi_2|^2
  \non
  \\
  & \quad {}
  + \cdelta \iO \an{\big|}\beta''(\phi_1)|\nabla\phi_1|^2 - \beta''(\phi_2)|\nabla\phi_2|^2\an{\big|}^2 \,.
  \non
\end{align}
\an{Besides}, we have that
\begin{align}
  & \iO |\beta'(\phi_1)\Delta\phi_1 - \beta'(\phi_2)\Delta\phi_2|^2 
  \leq c \iO |\phi|^2 \, |\Delta\phi_1|^2
  + c \iO |\Delta\phi|^2
  \non
  \\
  & \leq c \, \norma\phi_4^2 \, \norma{\Delta\phi_1}_4^2
  + c \iO |\Delta\phi|^2
  \leq c \, \normaV\phi^2 \, \norma{\phi_1}_{\Hx3}^2
  + c \iO |\Delta\phi|^2
  \non
  \\
  & \leq c \, \norma{\phi_1}_{\Hx3}^2 \iO |\nabla\phi|^2 
  + c \, \norma{\phi_1}_{\Hx3}^2 \, \norma\phi_{\L\infty H}^2
  + c \iO |\Delta\phi|^2\,,
  \non
\end{align}
and we observe at once that the function $t\mapsto\norma{\phi_1(t)}_{\Hx3}$ is estimated in $L^2(0,T)$ (see Remark~\ref{Piureg}).
This will be sufficient when we apply the Gronwall lemma to the first term of the last line after time integration,
and useful to estimate the time integral of the second one in terms of $u$ on account of \eqref{secondacontdep}.
In order to control the term containing $\beta''$\an{,} we observe that
\begin{align}
  & \beta''(\phi_1) |\nabla\phi_1|^2 - \beta''(\phi_2) |\nabla\phi_2|^2
  \non
  \\
  & = \bigl( \beta''(\phi_1) - \beta''(\phi_2) \bigr) |\nabla\phi_1|^2
  + \beta''(\phi_2) \nabla\phi \cdot \nabla\phi_1
  + \beta''(\phi_2) \nabla\phi_2 \cdot \nabla\phi\,,
  \non
\end{align}
so that
\begin{align}
  & \an{\big|}\beta''(\phi_1) |\nabla\phi_1|^2 - \beta''(\phi_2) |\nabla\phi_2|^2\an{\big|}^2
  \leq c \, |\phi|^2 \, |\nabla\phi_1|^4
  + c \, |\nabla\phi|^2 \, |\nabla\phi_1|^2
  + c \, |\nabla\phi_2|^2 \, |\nabla\phi|^2
  \non
  \\
  & \leq c \, \norma{\nabla\phi_1}_\infty^4 \, |\phi|^2
  + c \, \bigl( |\nabla\phi_1|^2 + |\nabla\phi_2|^2 \bigr) |\nabla\phi|^2 \,.
\non 
\end{align}
Therefore, we deduce that
\begin{align}
  & \iO \an{\big|}\beta''(\phi_1)|\nabla\phi_1|^2 - \beta''(\phi_2)|\nabla\phi_2|^2\an{\big|}^2
  \non
  \\
  & \leq c \, \norma{\nabla\phi_1}_\infty^4 \norma\phi_{\L\infty H}^2
  + c \iO \bigl( \norma{\nabla\phi_1}_\infty^2 + \norma{\nabla\phi_2}_\infty^2 \bigr) |\nabla\phi|^2\,, 
\non 
\end{align}
and we observe that \eqref{piuregbis} applied to $\phi_1$ ensures that,
after \an{integration over time}, the last two terms can be controlled,
the former by a direct estimate, and the latter via Gronwall's lemma.
\pier{The last three terms of \eqref{perterzacontdep} can be treated without any difficulty. Then,}
by integrating \eqref{perterzacontdep} with respect to time, 
combining all the inequalities we have  established, choosing $\delta$ small enough,
accounting for \eqref{primacontdep} and \eqref{secondacontdep}, and applying the Gronwall lemma, \pier{we arrive at}
\Beq
  \norma\phi_{\C0V\cap\L2{\Hx4}} \leq c \, \norma u_{\L2\Vp} \,.
  \label{terzacontdep}
\Eeq
This concludes the proof of the part of \eqref{contdep} concerning \an{the difference} $\an{\phi={}}\phi_1-\phi_2$.
It remains to prove the estimates regarding the other components of the solutions.
To this end, it suffices to write \pier{\eqref{terza} and \eqref{seconda}}
for both solutions, take the differences, and compare with the estimates already obtained.


\section{The control problem}
\label{CONTROL}
\setcounter{equation}{0}

In this section, we deal with the control problem presented in Section~\ref{STATEMENT}.
It is understood that in the entire section all of the assumptions made on the state system and the initial datum $\phiz$,
which ensure the validity of Theorems~\ref{Wellposedness} and~\ref{Contdep}, are \an{in} force,
as well as those related to the control problem, i.e., \last{conditions} \accorpa{hpalpha}{hpUad}.
Moreover, the generic constants denoted by $c$ can also depend on the latter.
\an{Here, in the control problem, the variable $u$ no longer plays just the role of a fixed and bounded source term,
but it enters \Pbl\ as a control variable, which justifies the introduction of the upper bound given by $M$.}
Furthermore, we recall that the cost functional, the space $\calU$, and the set $\Uad$ of admissible controls, 
are defined in \accorpa{cost}{Uad}.
We also~set
\Beq
  \calY := \H1\Wp \cap \L2W
  \label{defY}
\Eeq
and recall that, for every $u\in\calU$, there exists a unique solution $\soluz$ to problem \Pbl.
Since $\phi$ belongs to~$\calY$, we could define the {\it control-to-state\/} mapping as a mapping 
\last{from $\calU$ into~$\calY$ as \pier{$u\mapsto \phi$}}.
However, it is convenient to consider just a \nbh\ of the set $\Uad$ of admissible controls as the domain of~$\calS$.
So, we introduce the set $\UR$ and the \an{solution} map $\calS:\UR\to\calY$ by setting
\begin{align}
  & \UR := \graffe{u\in\calU:\ \norma u_\infty<R},
  \quad \hbox{where} \quad 
  R := \max\graffe{\norma\umin_\infty,\norma\umax_\infty} + 1,
  \label{defUR}
  \\
  & \UR \ni u \mapsto \calS(u) := \phi,
  \quad \hbox{where $\soluz$ is the \an{unique} solution}
  \non
  \\
  & \quad \hbox{to the problem \Pbl\ associated with $u$}.
  \label{defS}
\end{align}
We notice at once that we can apply Theorem~\ref{Wellposedness} with $M=R$  and Remark~\ref{Piureg}.
By also accounting for the regularity of~$\an{F}$ \an{in \eqref{hpBeta}}, we obtain, in particular,~that 
\begin{align}
  & \norma\phi_{\H1\Vp\cap\L\infty W\cap\L2{\Hx4}}
  + \norma\mu_{\L2V}
  + \norma w_{\L\infty H\cap\L2W}
  \non
  \\
  & \quad {}
  + \an{\somma i04 \norma{F^{(i)}(\phi)}_\infty }
  + \norma{\Delta f'(\phi)}_{\L\infty H}
  + \norma{\nabla f'(\phi)}_{\L2\Linfty}
  \leq c 
  \label{stimacoeffl}
\end{align}
for every $u\in\UR$.
This estimate will be  used throughout the whole section.
In the next three subsections, we prove the existence of an optimal control,
the \Frechet\ differentiability of \an{the solution mapping} $\calS$, and the \an{first-order} necessary condition for optimality
\an{presented} in Section~\ref{STATEMENT}.


\subsection{Existence of an optimal control}
\label{OPTIMUM}

\an{First, let us prove that the minimization problem introduced in \eqref{control} admits at least \juerg{one} solution.}
\Bthm
\label{Optimum}
There exists at least \juerg{one} optimal \last{control} for the control problem \eqref{control},
that is, there exists \juerg{some} $\ustar\in\Uad$ such~that
\Beq
  \calJ(\ustar,\calS(\ustar)) \leq \calJ(u,\calS(u))
  \quad \hbox{for every $u\in\Uad$}.
  \label{optimum}
\Eeq
\Ethm

\Bdim
We use the direct method \an{of the calculus of variations. In this direction, 
we observe that $\cal J$ is bounded from below as it is nonnegative. Then,} observing that $\Uad$ is nonempty,
we take a minimizing sequence $\graffe{\un}$ and, for every~$n$,
denote by $(\phin,\mun,\wn)$ the corresponding state.
Since $\Uad$ is bounded in~$\calU$, we \juerg{can assume without loss of generality} that
\Beq
  \un \to \ustar
  \quad \hbox{weakly star in~$\an{L^\infty(Q)}$}
  \non
\Eeq
for some \an{$\ustar\in\Uad$}, \an{as $\Uad$} is convex and closed\an{.}
Now, we notice that the stability estimate \eqref{stimacoeffl} holds true for $(\phin,\mun,\wn)$ \pier{and some constant independent of} $n$.
Thus, we \juerg{can assume without loss of generality that}, \an{as $n\to\infty$,}
\begin{alignat}{2}
  & \phin \to \phistar
  \quad &&\hbox{weakly star in $\H1\Vp\cap\L\infty W$},
  \non
  \\
  & \mun \to \mustar
  \quad &&\hbox{weakly in $\L2V$},
  \non
  \\
  & \wn \to \wstar
  \quad &&\hbox{weakly in $\L2W$},
  \non
\end{alignat}
for some \an{limit} triplet $\soluzstar$.
Since the stability estimate implies that $\graffe{\phin}$ is bounded in~$\LQ\infty$,
we can suppose that the nonlinearities involved in the state system are \Lip\ continuous.
On the other hand, the above convergence property of $\{\phin\}$, and well-known strong convergence results
(see, e.g., \cite[Sect.~8, Cor.~4]{Simon}), imply the strong convergence of~$\{\phin\}$ \an{to $\phistar$ as $n \to\infty$}, e.g., 
in $\LQ\infty$.
Thus, it is immediately seen that $\soluzstar$ is the solution to the state 
system corresponding to~$\ustar$\an{, that is,} $\phistar=\calS(\ustar)$.
Therefore, we also conclude~that
\Beq
  \lim_{n\to\infty} \calJ(\un,\phin)
  = \calJ(\ustar,\phistar).
  \non
\Eeq
\an{Now}, since \juerg{$\graffe\un$ is \an{supposed to be} a minimizing sequence},
the \rhs\ of th\an{e} equality is \an{a} minimum of~$\calJ$,
\last{meaning that} $\ustar$ is an optimal control.
\Edim


\subsection{The control-to-state operator}
\label{FRECHET}

This subsection is devoted to \an{analyze some differentiability properties} of the control-to-state mapping \an{$\cal S$} 
defined in~\eqref{defS}.
Namely, we \an{are going to show} its \Frechet\ differentiability \an{between suitable Banach spaces}.
\last{Prior to exploring this further, it is crucial to introduce two preliminary lemmas.}
\Blem
\label{Density}
The space $\H1W$ is dense in $\H1\Wp\cap\L2H$.
\Elem

\Bdim
Let $z\in\H1\Wp\cap\L2H$.
We \an{d}efine $\cj$ and $\zn$ by means of the eigenfunctions of the eigenvalue problem \EigenPbl\ as follows:
\Beq
  \cj(t) := (z(t),\ej)
  \quad \hbox{for $j\geq1$}
  \aand
  \zn(t) := \somma j1n \cj(t) \ej
  \quad \hbox{for $n\geq1$ and a.e. $t\in(0,T)$.}
  \non
\Eeq
Then $\cj\in H^1(0,T)$ for every $j$, and $\zn\in\H1W$ for every~$n$.
We claim that $\{\zn\}$ converges \last{strongly} to $z$ in $\H1\Wp\cap\L2H$.
\an{The strong convergence in $\L2 H$ is rather \last{straightforward}}\juerg{, because}
\Beq
  \lim_{n\to\infty} \norma{z-\zn}_{\L2H}^2
  = \lim_{n\to\infty} \sum_{j>n} \ioT |\cj(t)|^2 \, dt
  = 0 \,.
  \non
\Eeq
\pier{Note that} the function $\juerg{v\mapsto|\overline v|^2+\norma{\Delta v}^2}$ is the square of a norm in $W$ that is equivalent to the standard one, \pier{as one can see arguing by contradiction and using the compact embedding $W\emb V $ and the elliptic regularity theory.
Hence,}
the map $L:\juerg{v}\mapsto\graffe{(\juerg{v},\ej)}\an{_{j=1}^\infty}$ is a topological isomorphism from $W$ onto the
weighted space $\ell^2_\lambda$ defined~by
\begin{align}
  & \ell^2_\lambda := \Bigl\{ (\aj)_{j=1}^\infty:\ \somma j1\infty  |\pier{\lambdaj}\aj|^2 < + \infty \Bigr\},
  \quad \hbox{with the norm defined by} 
  \non
  \\
  & \norma{(\aj)}_\lambda^2 := |\an{a_1}|^2 + \somma j2\infty |\pier{\lambdaj}\aj|^2\,.
  \non
\end{align}
It follows that the adjoint operator $L^*:(\ell^2_\lambda)^*\to\Wp$ is a topological isomorphism. 
Therefore, every $\zeta\in\Wp$ \an{admits} a representation of the form
\Beq
  \zeta = \somma j1\infty \bj \ej
  \quad \hbox{with} \quad
  \somma j2\infty |\pier{\lambdaj^{-1}} \bj|^2 < + \infty \,.
  \non
\Eeq
Moreover, this representation is unique, and $\bj=\<\zeta,\ej>_W$ for every~$j$.
Furthermore, the function 
\Beq
  \Wp \ni \zeta
  \mapsto |\< \juerg{\zeta,e_1} >_W|^2 + \somma j2\infty  |\pier{\lambdaj^{-1}\,}\<\juerg{\zeta},\ej>_W|^2
  \non
\Eeq
is the square of a norm in $\Wp$ that is equivalent to the standard one.
\last{It thus} remains to prove that
\Beq
  \lim_{n\to\infty} \sum_{j>n} \pier{\lambdaj^{-2}} \ioT |\cj'(t)|^2 \, dt
  = 0\an{\,,}
  \non
\Eeq
\an{b}ut this is \an{readily} seen, since
\Beq
  \somma j2\infty \pier{\lambdaj^{-2}} \ioT |\cj'(t)|^2 \, dt
  \leq c \, \norma{\dt z}_{\L2\Wp}^2 \an{\leq c}\,.
  \non
\Eeq
\Edim

\Blem
\label{Bravogianni}
The identity
\Beq
  \< \dt\zeta(t) , \calN\zeta(t) >_W 
  = \frac 12 \, \ddt \, \normaVp{\zeta(t)}^2
  \quad \aat
  \label{bravogianni} 
\Eeq
holds true for every $\zeta\in\H1\Wp\cap\L2H$ satisfying $\zetabar=0$ \aet.
\Elem

\Bdim
Let $\zeta\in\H1\Wp\cap\L2H$ satisfy $\zetabar=0$ \aet.
By \an{Lemma \ref{Density}}, we can find a sequence $\graffe{\xi_n}$ in $\H1W$ that converges to $\zeta$ \an{strongly} in $\H1\Wp\cap\L2H$.
We \an{next} set $\zeta_n:=\xi_n-\overline\xi_n$, so that $\zetabar_n=0$. 
Th\an{us}, the sequence $\graffe{\zeta_n}$ converges to $\zeta-\zetabar=\zeta$ \an{strongly} in $\H1\Wp\cap\L2H$,
 \an{and} the identity \eqref{bravogianni} holds for~$\zeta_n$.
Thus, for every fixed $s\leq t$ in~$[0,T]$, we have~that
\Beq
  \int_s^t \< \dt\zeta_n(\tau) , \calN\zeta_n(\tau) >_W \, d\tau
  = \frac 12 \bigl(
    \normaVp{\zeta_n(t)}^2
    - \normaVp{\zeta_n(s)}^2
  \bigr).
  \non
\Eeq
On the other hand, since (the restriction of) $\calN$ is a continuous linear operator 
from the subspace of $H$ of the zero mean value functions into~$W$\an{,}
by the elliptic regularity theory,
the convergence of $\{\zeta_n\}$ to $\zeta$ in $\L2H$ implies that $\{\calN\zeta_n\}$ converges to $\calN\zeta$ \an{strongly} in $\L2W$.
Hence, we have that 
\Beq
  \lim_{n\to\infty} \int_s^t \< \dt\zeta_n(\tau) , \calN\zeta_n(\tau) >_W \, d\tau
  = \int_s^t \< \dt\zeta(\tau) , \calN\zeta(\tau) >_W \, d\tau \,.
  \non
\Eeq
Now, we observe that $V$ coincides with the \an{real} interpolation space $(W,H)_{1/2}$\last{,} the interpolation being understood in the sense of the trace method in Hilbert spaces.
We deduce that $\Vp=(H,\Wp)_{1/2}$, so that we have the continuous embedding
\Beq
  \H1\Wp \cap \L2H \emb \C0\Vp \,.
  \non
\Eeq
It follows that the convergence of $\{\zeta_n\}$ to $\zeta$ in $\H1\Wp\cap\L2H$ \an{as $n\to\infty$} implies the strong convergence in $\Vp$ of the pointwise values.
\last{As a result, it follows that}
\Beq
  \lim_{n\to\infty} \normaVp{\zeta_n(s)} = \normaVp{\zeta(s)}
  \aand
  \lim_{n\to\infty} \normaVp{\zeta_n(t)} = \normaVp{\zeta(t)}\,, 
  \non
\Eeq
and we conclude that
\Beq
  \int_s^t \< \dt\zeta(\tau) , \calN\zeta(\tau) >_W \, d\tau
  = \frac 12 \bigl(
    \normaVp{\zeta(t)}^2
    - \normaVp{\zeta(s)}^2
  \bigr).
  \non
\Eeq
Since $s$ and $t$ are arbitrary in~$[0,T]$, this is equivalent to \eqref{bravogianni}, and the proof is complete.
\Edim

\an{
\Brem
\label{Bravogiannibis}
Let us remark that, arguing along the same lines as above, \pier{\last{the} analogous formula}
\Beq
  \< \dt z , \Delta z >_W
  = - \frac 12 \, \ddt \, \norma{\nabla z}^2 \juerg{\quad\mbox{a.e. in $(0,T)$}},
  \label{bravogiannibis}
\Eeq
\pier{can be shown} for every  $z\in\H1\Wp\cap\L2W$ such that $\Delta z\in\L2W$.
Let us sketch the corresponding proof. \last{First, we set}
\Beq
  W_1 := \graffe{v\in W: \ \Delta v\in W},
  \non
\Eeq
and consider the formal sum
\Beq
  \somma j1\infty \cj \ej \,.
  \label{sum}
\Eeq
Then, \eqref{sum} represents an element of $W$ if and only if
$\somma j1\infty|\lambdaj\cj|^2<+\infty$.
\last{Besides}, it follows that \eqref{sum} characterizes an element of either $W_1$, or $\Wp$, or $V$,
if and only~if
\Beq
  \somma j1\infty |\lambdaj^2 \cj|^2 < +\infty\,,
  \quad \hbox{or} \quad
  \somma j2\infty |\lambdaj^{-1} \cj|^2 < +\infty\,,
  \quad \hbox{or} \quad
  \somma j1\infty |\lambdaj^{1/2} \cj|^2 < +\infty\,,
  \non
\Eeq
respectively. \juerg{It follows} that the interpolation space $(W_1,\Wp)_{1/2}$ coincides with~$V$,
implying the continuity of the embedding
\Beq
  \H1\Wp \cap \L2{W_1} \emb \C0V \,.
  \non
\Eeq
\last{Finally}, the same argument \juerg{as that} given in the proof of the above lemma
entails \eqref{bravogiannibis}.
\Erem
}

Our result on the \Frechet\ differentiability of $\calS$ involves the linearized system we introduce \last{in the subsequent lines}.
For fixed $u\in\UR$ and $h\in \an{\L2 H}$, the linearized system corresponding to $u$ and the variation~$h$
consists in looking for a triplet $\soluzl$ satisfying
\begin{align}
  & \soluzl \in \bigl( \H1\Wp\cap\L2W \bigr) \times \L2H \times \L2W\,,
  \label{regsoluzl}
  \\
  & \< \dt\psi , v >_W 
  \an{
  - \iO \eta \, \Delta v
  + \sigma \iO \psi v 
  }
  = \iO h v
  \quad \hbox{\aet, for every $v\in W$}\,,
  \label{primal}
  \\
  & \iO \nabla\omega \cdot \nabla v
  + \iO f''(\phi) \, w  \psi  v
  + \iO f'(\phi) \, \omega  v
  + \nu \iO \omega v
  = \iO \eta  v
  \non
  \\
  & \quad \hbox{\aet, for every $v\in V$}\,,
  \label{secondal}
  \\
  & \iO \nabla\psi \cdot \nabla v + \iO f'(\phi) \, \psi v
  = \iO \omega v
  \quad \hbox{\aet, for every $v\in V$}\,,
  \label{terzal}
  \\
  & \psi(0) = 0\,, 
  \label{cauchyl}
\end{align}
\Accorpa\Pbll primal cauchyl
where $\phi$ and $w$ are the components of the solution $\soluz$ to the state system associated with~$u$.
We observe that the continuous embedding $\H1\Wp\cap\L2W\emb\C0H$ 
gives a precise meaning to the initial condition~\eqref{cauchyl}.
We also notice that \eqref{secondal} and \eqref{terzal} could have been written in the form
\Beq
  - \Delta\omega + f''(\phi) \, w \, \psi + f'(\phi) \, \omega \,\juerg{+\nu\,\omega}
  = \eta 
  \aand
  - \Delta\psi + f'(\phi) \, \psi 
  = \omega
  \quad \aeQ\,,
  \non
\Eeq
due to the regularity of $\omega$ and $\psi$ given by \eqref{regsoluzl}, \an{because the space $W$ \last{inherently} \juerg{encodes} the boundary conditions}.
\an{Nevertheless, in \juerg{the following} we prefer} to keep \juerg{on} \an{working with} the \an{aforementioned} variational formulation.

\Bthm
\label{Wellposednessl}
For every $u\in\UR$ and $h\in\an{\L2 H}$, problem \Pbll\ has a unique solution $\soluzl$ \an{with the regularity specified 
in \eqref{regsoluzl}}.
Moreover, the estimate
\Beq
  \norma\psi_{\H1\Wp\cap\L2W}
  + \norma\eta_{\L2H}
  + \norma\omega_{\L2W}
  \leq K_3 \, \norma h_{\L2H}
  \label{continuita}
\Eeq
\pier{holds true with a positive constant $K_3$ that depends only on $\Omega$, $T$, the structure of the system,
and \an{the upper bound $R$.}}
\Ethm

\Bdim
We recall \eqref{stimacoeffl}, and we start by proving uniqueness.
\last{It is worth noting that \juerg{a portion} of the insights we are about to unveil also contributes significantly to the proof of existence.}

\step
First a priori estimate

Testing \eqref{primal} by $1/|\Omega|$ yields an ordinary differential equation for the mean value~$\psibar$, namely
\Beq
  \ddt \, \psibar + \sigma \, \psibar = \hbar \juerg{\quad\mbox{a.e. in $(0,T)$}}\,.
  \label{primalbar}
\Eeq
\an{From \eqref{cauchyl}, we \juerg{infer that} $\psibar(0)=0$. Thus,} we \last{readily} deduce that
\Beq
  \norma\psibar_{\pier{H^{1}(0,T)}} \leq c \, \norma h_{\L2H} \,.
  \label{primastimal}
\Eeq

\step
Second a priori estimate

We \juerg{first} multiply \eqref{primalbar} by a generic $v\in W$ and integrate over~$\Omega$.
Then, we take the difference between \eqref{primal} and the equality just obtained
and test it by~$\calN(\psi-\psibar)$ (recall~\eqref{defN}).
Thanks to Lemma~\ref{Bravogianni}, \an{along with} the definition of~$\calN$, we \an{find}~that
\Beq
  \frac 12 \, \ddt \, \normaVp{\psi-\psibar}^2
  \an{{}
  + \iO \eta \, (\psi-\psibar)
  + \sigma \, \normaVp{\psi-\psibar}^2
  {}}
  = \iO h \, \calN(\psi-\psibar) \juerg{\quad\mbox{a.e. in $(0,T)$}}\,.
  \non
\Eeq
At the same time, we test \eqref{secondal} and \eqref{terzal} by $\psi-\psibar$ and $\psi-\omega$, respectively,
and rearrange \last{the terms} \an{to} obtain\an{, in the order,}~that
\begin{align}
  & \iO \nabla\omega \cdot \nabla\psi 
  = \iO \eta(\psi-\psibar)
  - \iO \bigl( f''(\phi) \, \psi \, w + f'(\phi) \, \omega + \nu \omega \bigr) (\psi-\psibar)\,,
  \non
  \\
  & \iO |\nabla\psi|^2
  - \iO \nabla\psi \cdot \nabla\omega
  + \iO |\omega|^2
  = \iO \omega \, \psi
  - \iO f'(\phi) \, \psi \, (\psi-\omega) \,.
  \non
\end{align}
At this point, we add the above equalities to each other, \an{noticing that some cancellations occur\pier{. We} infer}~that
\begin{align}
  & \frac 12 \, \ddt \, \normaVp{\psi-\psibar}^2
  + \sigma \, \normaVp{\psi-\psibar}^2
  + \iO |\nabla\psi|^2
  + \iO |\omega|^2
  \non
  \\
  & = \iO h \, \calN(\psi-\psibar) 
  - \iO f''(\phi) \, \psi \, w \, (\psi-\psibar)
  + \iO \omega \, \psi
  \non
  \\
  & \quad {}
  + \iO f'(\phi) \, \omega \, \psibar 
  - \iO f'(\phi) \, |\psi|^2 
  \pier{{}+ \nu \iO \omega (\psi- \psibar)} \juerg{\quad\mbox{a.e. in $(0,T)$}}\,.
  \label{persecondastimal}
\end{align}
The \rhs\ is easily dealt with, in particular, by using the compacness inequality \eqref{compact} 
to handle the $H$ norm of~\pier{$\psi-\psibar$}, since one can owe to the Gronwall lemma after time integration.
We just \an{detail how to handle} the most delicate term, i.e., \an{the one involving}~$w$.
On account of \eqref{primastima}, it suffices to estimate the analogue obtained by replacing $\psi$ by $\psi-\psibar$.
We have~that
\Beq
  - \iO f''(\phi) \, w \, |\psi-\psibar|^2
  \leq \norma{f''(\phi)}_\infty \, \norma w \, \norma{\psi-\psibar}_4^2 \,.
  \non
\Eeq
Observing that the product $\norma{f''(\phi)}_\infty\,\norma w$ is estimated in $L^\infty(0,T)$ by virtue of~\eqref{stimacoeffl},
and owing to the compactness inequality \eqref{compact}, we deduce~that
\Beq
  - \iO f''(\phi) \, w \, |\psi-\psibar|^2
  \leq \frac 12 \iO |\nabla\psi|^2
  + c \, \normaVp{\psi-\psibar}^2 \juerg{\quad\mbox{a.e. in $(0,T)$}}\,.
  \non
\Eeq
Therefore, after integrating \eqref{persecondastimal} with respect to time,
we apply Gronwall's lemma as announced and obtain\an{,} recalling \an{also} \eqref{primastimal}, \an{that}
\Beq
  \norma\psi_{\L\infty\Vp\cap\L2V}
  + \norma\omega_{\L2H}
  \leq c \, \norma h_{\L2H} \,.
  \label{secondastimal}  
\Eeq

\step
Consequence

By applying \eqref{secondastimal} in the case $h=0$, we obtain that $\psi=\omega=0$.
Then \eqref{secondal} yields that $\eta=0$.
Since \an{the system} \Pbll\ is linear, this proves the uniqueness part of the statement.

It remains to prove the existence of a solution and the estimate~\eqref{continuita}.
This can be done by a discretization procedure \an{like the} Faedo--Galerkin scheme, 
with the same set of eigenfunctions as that used for Theorem~\ref{Wellposedness}:
the unknowns of the discrete problem are $\Vn$-valued, 
and the analogues of \accorpa{primal}{terzal} are satisfied for every $v\in\Vn$.
One obtains a linear system that can be solved on account of the \an{zero} initial condition derived from~\eqref{cauchyl}.
Then, one performes a priori estimates \pier{and let $n$ tend} to infinity.
\last{However, for the sake of brevity, we limit ourselves to demonstrating formal estimates, acknowledging that the \an{test} functions we employ \juerg{are} admissible at the discrete level, especially when only the variational form of the problem is considered since} $\Delta^k v\in\Vn\subset W$ for every $k$, $n$ and $v\in\Vn$.
Since also estimates \eqref{primastimal} and~\eqref{secondastimal} can be performed at the discrete level, 
we can account for them in the sequel.

\step
Third a priori estimate

We test \eqref{terzal} by $-\Delta\psi$.
By virtue of \eqref{secondastimal}\an{,} we derive an estimate for $\Delta\psi$ in $\L2H$;
\juerg{from elliptic regularity theory we then infer that}
\Beq
  \norma\psi_{\L2W} \leq c \, \norma h_{\L2H} \,.
  \label{terzastimal}
\Eeq

\step
Fourth a priori estimate

We test \eqref{primal}, \eqref{secondal}, and \eqref{terzal}, by $\psi$, $-\Delta\psi$, and~$\Delta\omega$, respectively.
Then, we sum up and rearrange.
Thanks to the occurring cancellations and an integration by parts in the terms involving~$\Delta\omega$, we obtain~that
\begin{align}
  & \frac 12 \, \ddt \iO |\psi|^2
  + \sigma \iO |\psi|^2 
  + \iO |\nabla\omega|^2
  \non
  \\
  & = \iO h \psi
  \pier{{}-{}} \iO \bigl( f''(\phi) \, w  \psi  + f'(\phi) \, \omega + \nu \omega \bigr) (-\Delta\psi)
  - \iO \bigl( f''(\phi) \, \psi \nabla \phi + f'(\phi) \nabla\psi \bigr) \cdot \nabla\omega\,,
  \non
\end{align}
and we \an{are left with} estimat\an{ing} the \rhs.
We just consider the most delicate terms, since the others are easier \juerg{and can be left as an  
exercise to the interested reader}.
We recall \eqref{stimacoeffl} and \eqref{terzastimal} and \an{point out}~that
\begin{align}
  & \pier{{}-{}} \iO f''(\phi) \, w  \psi  (-\Delta\psi)
  - \iO f''(\phi) \, \psi \nabla \phi \cdot \nabla\omega
  \non
  \\
  & \leq c \, \norma w^2 \, \norma\psi_\infty^2
  + c \, \norma{\Delta\psi}^2
  + \frac 12 \iO |\nabla\omega|^2
  + c \, \norma{\nabla\phi}^2 \, \norma\psi_\infty^2
  \non
  \\
  & \leq c \, \normaW\psi^2
  + \frac 12 \iO |\nabla\omega|^2 \,.
  \non
\end{align}
Hence, after integrating \an{over} time, and also \an{accounting for} \eqref{secondastimal} and \eqref{terzastimal},
we conclude that
\Beq
  \norma\psi_{\L\infty H} + \norma\omega_{\L2V}
  \leq c \, \norma h_{\L2H} \,.
  \label{quartastimal}
\Eeq

\step
Fifth a priori estimate

We test \eqref{primal}, \eqref{secondal}, and \eqref{terzal}, by
$-\Delta\psi$, $\Delta^2\psi$, and $-\Delta^3\psi$, respectively,
sum up, perform some integrations by parts, and account for some cancellations.
\an{This leads us to the identity}
\begin{align}
  & \frac 12 \, \ddt \iO |\nabla\psi|^2
  + \sigma \iO |\nabla\psi|^2
  + \iO |\Delta^2\psi|^2
  \non
  \\
  & {} = 
  - \iO h \Delta\psi
  \pier{{}-{} \nu \iO \omega \Delta^2\psi}
  - \iO f''(\phi) \, w \psi \Delta^2\psi
    \non
  \\
  & \quad {}
  - \iO f'(\phi) \, \omega \Delta^2\psi
  + \iO f'(\phi) \, \psi \Delta^3\psi \,.
  \label{perquintastimal}
\end{align}
As before, we just handle the terms on the \rhs\ that need some treatment\last{, while leaving the simpler ones for the reader to follow.}
Integrating with respect to time, and owing to \eqref{stimacoeffl} and \eqref{quartastimal}, we obtain~that
\begin{align}
  & - \intQt f''(\phi) \, w \psi \Delta^2\psi
  \leq c \iot \norma{w(s)}_\infty \, \norma{\psi(s)} \, \norma{\Delta^2\psi(s)} \, ds
  \non
  \\
  & \leq \frac 14 \intQt |\Delta^2\psi|^2
  + c \, \norma w_{\L2W}^2 \, \norma\psi_{\L\infty H}^2
  \leq \frac 14 \intQt |\Delta^2\psi|^2
  + c \, \norma h_{\L2H}^2 \,.
  \non
\end{align}
Next, \last{integrating with respect to time,} \juerg{invoking} \eqref{stimacoeffl} once more and the above estimates, we have~that
\begin{align}
  & \last{\intQt} f'(\phi) \, \psi \Delta^3\psi
  = \last{\intQt}\bigl( (\Delta f'(\phi)) \psi + 2 (\nabla f'(\phi)) \cdot \nabla\psi + f'(\phi) \Delta\psi \bigr) \Delta^2\psi
  \non
  \\
  & \leq \frac 14 \intQt |\Delta^2\psi|^2
  + c \, \norma{\Delta f'(\phi)}_{\L\infty H}^2 \, \norma\psi_{\L2W}^2
  \non
  \\
  & \quad {}
  + c \iot \norma{\nabla f'(\phi(s)\an{)}}_\infty^2 \, \norma{\nabla\psi(s)}^2 \, ds
  + c \intQt |\Delta\psi|^2
  \non
  \\
  & \leq \frac 14 \intQt |\Delta^2\psi|^2
  + c \, \norma h_{\L2H}^2
  + c \iot \norma{\nabla f'(\phi(s)\an{)}}_\infty^2 \, \norma{\nabla\psi(s)}^2 \, ds\,,
  \non
\end{align}
\an{and} we observe that the function $t\mapsto\norma{\nabla f'(\phi(s))}_\infty^2$ 
is estimated in $L^1(0,T)$.
Therefore, \pier{after integration of \eqref{perquintastimal} and application of} Gronwall's lemma, we conclude~that
\Beq
  \norma{\nabla\psi}_{\L\infty H}
  + \norma{\Delta^2\psi}_{\L2H}
  \leq c \, \norma h_{\L2H} \,,
  \non
\Eeq
whence\an{, by elliptic regularity theory,}
\Beq
  \norma\psi_{\juerg{L^\infty(0,T;V)}\cap\L2{\Hx4}}
  \leq c \, \norma h_{\L2H} \,.
  \label{quintastimal}
\Eeq

\step
Consequences

By testing \eqref{terzal} by $\Delta^2\omega$ and integrating by parts we get the identity
\juerg{
\begin{align}
\iO |\Delta\omega|^2\,=\,-\iO \Delta^2\psi\,\Delta\omega\,+\iO\Delta(f'(\phi)\psi)\,\Delta\omega\,,
\non
\end{align}
whence, integrating over time and invoking Young's inequality,
\begin{align*}
\intQt |\Delta\omega|^2\,\le\,c\intQt |\Delta^2\psi|^2\,+\,c\intQt |\Delta(f'(\phi)\psi)|^2\,.
\end{align*}
Using \eqref{stimacoeffl} and the estimates proved above, it is not difficult to see that the \rhs\ is bounded by
$\,c\,\|h\|^2_{L^2(0,T;H)}$, whence, by elliptic regularity, }
\Beq
  \norma\omega_{\L2W} \leq c \, \norma h_{\L2H} \,.
  \label{stimaomega}
\Eeq
Testing \eqref{secondal} by $\eta$ \juerg{then easily} yields~that
\Beq
  \norma\eta_{\L2H} \leq c \, \norma h_{\L2H} \,.
  \label{stimaeta}
\Eeq
Finally, by comparison in \eqref{primal}
(at~the discrete level one should \an{argue} as we did to obtain~\eqref{sestastima}),
we~readily find~that
\Beq
  \norma{\dt\psi}_{\L2\Wp}
  \leq c \, \norma h_{\L2H} \,.
  \label{stimadtpsi}
\Eeq
This concludes the proof.
\Edim

We are now ready to prove the \Frechet\ differentiability of  the map $\calS$ defined in~\eqref{defS}.

\Bthm
\label{Frechet}
The control-to-state mapping $\calS$ is \Frechet\ differentiable in $\UR$
as a mapping from \juerg{$\UR\subset\calU$} into $\calY$.
Moreover, given $u\in\UR$, the \Frechet\ derivative $D\calS(u)\juerg{\in \calL(\calU,\calY)}$ is the linear operator
that maps any $h\in\calU$ into the component $\psi$ of the solution $\soluzl$ to the linearized problem \Pbll\
associated with $u$ and the variation~$h$.
\Ethm

\Bdim
The map described in the statement actually belongs to $\calL(\calU,\calY)$ as a consequence of~\eqref{continuita} 
\an{in Theorem \ref{Wellposednessl}. To prove the claim, we} fix some $u\in\UR$ and \juerg{show that}
\Beq
  \lim_{\norma h_\infty\to0} \frac {\norma{\calS(u+h)-\calS(u)-\psi}_\calY} {\norma h_\infty}
  = 0\,,  \label{frechet} 
\Eeq
\an{where $\psi$ is the first component of the solution to the linearized associated to $u$ and $h$.}
Without loss of generality, we assume that $\norma h_\infty<R-\norma u_\infty$,
so that $u+h$ also belongs to \an{the open set}~$\UR$.
We denote by $\soluz$ and $(\phih,\muh,\wh)$ the solutions to the state system corresponding to $u$ and~$u+h$, respectively,
and notice that the uniform bound \eqref{stability} given by Theorem~\ref{Wellposedness} holds true for both of them.
Moreover, $\soluzl$ \an{denotes} the solution to the linearized system as in the statement, \an{and}
we set for convenience
\Beq
  \psih := \phih - \phi - \psi \,, \quad
  \etah := \muh - \mu - \eta\,,
  \aand
  \omegah := \wh - w - \omega \,.
  \non
\Eeq
We aim at proving an inequality that implies~\eqref{frechet}.
Namely\an{, we are going to prove that}
\Beq
  \norma\psih_{\H1\Wp\cap\L2W} + \norma\omegah_{\L2V}
  \leq c \, \norma h_{\L2\Vp}^2 \,,
  \label{tesi}
\Eeq
\an{where the first norm, due to the above definitions, is nothing but the numerator of \eqref{frechet}.}
To this end, we observe that $(\psih,\etah,\omegah)$ solves the following problem:
\begin{align}
  & \< \dt\psih , v >_W 
  \an{{}
  - \iO \etah \, \Delta v
  + \sigma \iO \psih v
  }
  = 0
  \quad \hbox{\aet, for every $v\in W$}\,,
  \label{primah}
  \\
  & \iO \nabla\omegah \cdot \nabla v
  + \nu \iO \omegah v
  = \iO \etah  v
  - \iO \Lambda_1 \, v
  \non
  \\
  & \quad \hbox{\aet, for every $v\in V$}\,,
  \label{secondah}
  \\
  & \iO \nabla\psih \cdot \nabla v
  = \iO \omegah v
  - \iO \Lambda_2 \, v
  \quad \hbox{\aet, for every $v\in V$}\,,
  \label{terzah}
  \\
  & \psih(0) = 0 \,,
  \label{cauchyh}
\end{align}
where we have set for brevity
\begin{align}
  \Lambda_1 & :=
  f'(\phih) \wh - f'(\phi) w - f''(\phi) \psi w - f'(\phi) \omega
  \non
  \\
  &  {\phantom{:}} = [ f'(\phih) - f'(\phi) - f''(\phi) \psi ] w
  + [ f'(\phih) - f'(\phi) ] (\wh-w)
  + f'(\phi) \omegah
  \non
  \\
  \Lambda_2 & := 
  f(\phih) - f(\phi) - f'(\phi) \psi \,.
  \non
\end{align}
We notice at once that expanding $f'$ and $f$ around~$\phi$
(namely, around $\phi(x,t)$ for almost all fixed $(x,t)$ in~$Q$)
by means of Taylor's formula with \pier{integral} remainder, \an{we find} that
\begin{align}
  & f'(\phih) - f'(\phi) - f''(\phi) \psi
  = f''(\phi) \psih + R_1 \,,
  \non
  \\
  & f(\phih) - f(\phi) - f'(\phi) \psi
  = f'(\phi) \psih + R_2 \,,
  \non
\end{align}
where \an{the remainders} $R_1$ and $R_2$ have the form\pier{%
\begin{align}
 & R_1 =\Bigl[ \int_0^1 (1-s) f'''(s\phih+(1-s)\phi)) \, ds \Bigr] (\phih-\phi)^2
  \, ,  \non
  \\
  & R_2 = \Bigl[ \int_0^1 (1-s) f''(s\phih+(1-s)\phi)) \, ds \Bigr] (\phih-\phi)^2 \, .
  \non
\end{align}}%
Hence, we have that
\begin{align}
  & \Lambda_1 
  = [ f''(\phi) \psih + R_1 ] w
  + [ f'(\phih) - f'(\phi) ] (\wh-w)
  + f'(\phi) \omegah\,,
  \non
  \\
  & \aand
  \Lambda_2 
  = f'(\phi) \psih + R_2 \,,
  \label{formuleLambda}
\end{align}
where the remainders satisfy
\Beq
  |R_1| \leq c \, (\phih-\phi)^2
  \aand
  |R_2| \leq c \, (\phih-\phi)^2
  \quad \aeQ \,,
  \label{stimaR}
\Eeq
since both $\phi$ and $\phih$ are uniformly bounded in $L^\infty$ as a consequence of~\eqref{stability}
and \pier{the integration variable $s$ \last{attains} values} in~$[0,1]$.
\an{In the estimates to follow, it} is understood that $\delta$ is a positive parameter \an{whose value will be chosen at the end of the computations}.

\step
First a priori estimate

\juerg{First observe} that $\psih$ has zero mean value, \an{which directly follows from testing \eqref{primah} 
by $1/|\Omega|$ and using~\eqref{cauchyh}}. \juerg{Next,}
we test \eqref{primah}, \eqref{secondah}, and \eqref{terzah}, 
by $\calN\psih$, $\psih$, and $\psih-\omegah$, respectively, 
and add the \juerg{resulting equalities} to each other.
By noting some cancellations, and recalling Lemma~\ref{Bravogianni}, we have that, \juerg{a.e.~in~$(0,T)$}, 
\begin{align}
  & \frac 12 \ddt \, \normaVp\psih^2
  + \sigma \, \normaVp\psih^2
  + \iO |\nabla\psih|^2
  + \iO |\omegah|^2
  \non
  \\
  & = - \iO \Lambda_1 \, \psih
  + \iO \Lambda_2 \, (\omegah-\psih)
  + (1-\nu) \iO \omegah \psih \,.
  \label{perprimaF}
\end{align}
\an{As the \juerg{last three} terms on the \lhs\ are nonnegative, we just} have to estimate the \an{integrals on the} \rhs.
By the expression of $\Lambda_1$ given in \eqref{formuleLambda}, 
\an{the term $\iO \Lambda_1 \, \psih$ may be written as} sum of three contributions\juerg{, namely,}
\begin{align}
  & I_1 := \iO [ f''(\phi) \psih + R_1 ] \, w \psih \,, \quad
  I_2 := \iO [ f'(\phih) - f'(\phi) ] (\wh-w) \psih\,,
  \non
  \\& 
  \aand
  I_3 := \iO f'(\phi) \omegah \psih \,.
  \non
\end{align}
\juerg{Now, observe that} \pier{$I_3$}, as well as the last term of \eqref{perprimaF}, 
can easily be treated using the Young and compactness inequalities.
The latter inequality is also used to \pier{deal with}~$I_1$.
\pier{Indeed, on} account of \eqref{stimaR}, we have~that
\begin{align}
  & |I_1|
  \leq c \, I_1' + c \, I_1''\,,
  \quad \hbox{where} \quad
  I_1' := \iO |w| \, |\psih|^2
  \aand
  I_1'' := \iO (\phih-\phi)^2 \, |w| \, |\psih|\,. 
  \non
\end{align}
\juerg{The term} $I_1'$ can immediately be estimated \pier{with the help of  \eqref{stability} and \eqref{piureg} (cf.~Remark~\ref{Piureg})} \an{as}
\Beq
  I_1'
  \leq \norma w_{\L\infty H} \, \norma\psih_4^2
  \leq \delta \, \iO |\nabla\psih|^2 + \cdelta \, \normaVp\psih^2\,, 
  \non
\Eeq
while the estimation of $I_1''$ \an{requires} more work.
We estimate its time integral by first \pier{using} the \Holder, Sobolev, and Poincar\'e inequalities,  \juerg{where we recall} that $\psih$ has zero mean value.
Then, we apply the stability estimate for $w$ related to \eqref{piureg}, Young's inequality, 
and the continuous dependence inequality \eqref{contdep} with $u_1=u+h$ and $u_2=u$.
\pier{Thus, we find out}~that
\begin{align}
  & \intQt |\phih-\phi|^2 \, |w| \, |\psih| 
  \leq \iot \norma{\phih(s)-\phi(s)}_6^2 \, \norma{w(s)}_2 \, \norma{\psih(s)}_6 \, ds
  \non
  \\
  & \leq c \iot \normaV{\phih(s)-\phi(s)}^2 \, \norma{w(s)} \, \normaV{\psih(s)} \, ds
  \non
  \\
  & \leq c \iot \normaV{\phih(s)-\phi(s)}^2 \, \norma{w(s)} \, \norma{\nabla\psih(s)} \, ds
  \non
  \\
  &  \leq \frac 12 \intQt |\nabla\psih|^2
  + c \, \norma w_{\L\infty H}^2 \, \norma{\phih-\phi}_{\L\infty V}^4
  \non
  \\
  & \leq \frac 12 \intQt |\nabla\psih|^2
  + c \, \norma h_{\L2\Vp}^4 \,\an{,}
  \label{bravoandrea}
\end{align}
\an{completing} the estimate of~$|I_1|$.
Next, we observe that
\Beq
  |I_2|
  \leq \norma\psih^2
  + c \, \norma{\phih-\phi}_4^2 \, \norma{\wh-w}_4^2
  \non
\Eeq
and we can \juerg{apply} the compactness inequality once more to the first term,
while we observe that \eqref{contdep} ensures~that
\begin{align}
  & \last{\iot} \norma{\phih(\last{s})-\phi(\last{s})}_4^2 \, \norma{\wh(\last{s})-w(\last{s})}_4^2 \, \last{ds}
  \non
  \\
  & \leq c \, \norma{\phih-\phi}_{\L\infty V}^2 \, \norma{\wh-w}_{\L2V}^2
  \leq c \, \norma h_{\L2\Vp}^4 \,.
  \non
\end{align}
Now, we consider the term involving~$\Lambda_2$.
By recalling \eqref{formuleLambda} and \eqref{stimaR}, we have~that
\begin{align}
  & \iO \Lambda_2 \, (\omegah-\psih)
  = \iO \bigl( f'(\phi) \psih + R_2 \bigr) (\omegah-\psih)
  \leq c \iO \bigl( |\psih| + (\phih-\phi)^2 \bigr) (|\omegah|+|\psih|).
  \non
\end{align}
Hence, it can be estimated with the same ideas as we used to prove \eqref{bravoandrea} and is even easier.
At this point, we return to \eqref{perprimaF}, integrating it with respect to time.
By choosing $\delta$ small enough, and applying the Gronwall lemma, we conclude~that
\Beq
  \norma\psih_{\L\infty\Vp\cap\L2V} 
  + \norma\omegah_{\L2H}
  \leq c \, \norma h_{\L2\Vp}^2 \,.
  \label{primaF}
\Eeq
\juerg{Then,} \pier{using \eqref{primaF} to estimate the \rhs\ of}  \eqref{terzah}, we deduce~that
\Beq
  \norma\psih_{\L2W}
  \leq c \, \norma h_{\L2\Vp}^2 \,.
  \label{daprimaF}
\Eeq

\step
Second a priori estimate

\juerg{For the sake of brevity, we argue only formally}.
However, we notice that $(\psih,\etah,\omegah)$
can be approximated by $\Vn$-valued functions (recall~\eqref{defVn}),
since this is true for $(\phih,\muh,\wh)$, $\soluz$, and $\soluzl$.
Hence, the estimate we are going to \an{perform}
can be \juerg{justified rigorously.}
We test \eqref{primah}, \eqref{secondah}, and \eqref{terzah},
by $\psih$, $-\Delta\psih$, and~$\Delta\omegah$,
add the resulting identities, and account for some cancellations.
We then obtain~that \juerg{a.e. in $(0,T)$ it holds}
\Beq
  \frac 12 \ddt \, \norma\psih^2
  + \sigma \, \norma\psih^2
  + \norma{\nabla\omegah}^2
  = \nu \iO \omegah \Delta\psih
  + \iO \Lambda_1 \, \Delta\psih
  - \iO \Lambda_2 \, \Delta\omegah \,.
  \label{persecondaF}
\Eeq
The first term on the \rhs\ can easily be treated 
by the \Holder\ inequality and~\accorpa{primaF}{daprimaF}.
As for the next term, we apply the Young inequality and \eqref{daprimaF}\juerg{.
We then} see that we are led to treat the integral of \pier{$ |\Lambda_1|^2$}.
By recalling \eqref{formuleLambda} and~\eqref{stimaR}, we have~that
\Beq
  \iO |\Lambda_1|^2
  \leq c \, \norma w_\infty^2 \, \norma\psih^2 
  + c \, \norma{\phih-\phi}_4^4 \, \norma w_\infty^2
  + c \, \norma{\phih-\phi}_4^2 \, \norma{\wh-w}_4^2 
  + c \, \norma\omegah^2 \,.
  \non
\Eeq
By noting that the function $s\mapsto\norma{w(s)}_\infty^2$ is bounded in $L^1(0,T)$ by \eqref{stability},
we can treat the first addend by the Gronwall lemma after time integration.
On the other hand, the time integrals of the remaining terms on the \rhs\ 
can be treated by owing to  \eqref{stability}, \eqref{contdep} and~\eqref{primaF}.
\an{Namely, we} conclude~that
\Beq
  \last{\intQt} |\Lambda_1|^2
  \leq c \, \norma h_{\L2\Vp}^4 \,.
  \non
\Eeq
Finally, by \eqref{formuleLambda} and~\eqref{stimaR}, we have~that
\Beq
  - \iO \Lambda_2 \, \Delta\omegah
  \leq \delta \, \norma{\nabla\omegah}^2
  + \cdelta \, \norma{\nabla\Lambda_2}^2
  \leq \delta \, \norma{\nabla\omegah}^2
  + \cdelta \, \pier{\norma{\psih}_V^2}
  + \cdelta \, \norma{\nabla R_2}^2 \,.
  \non
\Eeq
Since we can owe to \eqref{daprimaF}, it remains to estimate the last term.
\pier{To this aim, we compute $\nabla R_2$ by differentiating, under the integral sign as well, and \last{subsequently obtain}}
\begin{align}
  & \nabla R_2
  = 2 (\phih-\phi) \nabla(\phih-\phi) \int_0^1 (1-s) f''(s\phih+(1-s)\phi)) \, ds 
  \non
  \\
  & \quad {}
  + (\phih-\phi)^2 \int_0^1 (1-s) f'''(s\phih+(1-s)\phi)) \bigl( s\nabla\phih + (1-s)\nabla\phi \bigr) \, ds\,\pier{.}
  \non
\end{align}
\pier{Then,} we deduce~that
\Beq
  |\nabla R_2|
  \leq c \, |\phih-\phi| \, |\nabla(\phih-\phi)|
  + c \, \bigl( |\nabla\phih|+|\nabla\phi| \bigr) (\phih-\phi)^2
  \quad \aeQ .
  \non
\Eeq
Therefore, \pier{it turns out that}
\Beq
  \norma{\nabla R_2}^2
  \leq c \, \norma{\phih-\phi}_4^2 \, \norma{\nabla(\phih-\phi)}_4^2
  + c \, \bigl( \norma{\nabla\phih}_6^2 + \norma{\nabla\phi}_6^2 \bigr) \, \norma{\phih-\phi}_6^4\,,
  \non
\Eeq
\pier{from which we infer}
\begin{align}
  & \last{\iot} \norma{\nabla R_2(s)}^2 \, ds
  \leq c \, \norma{\phih-\phi}_{\C0V}^2 \, \norma{\phih-\phi}_{\L2W}^2
  \non
  \\
  & \quad {}
  + c \, \bigl( \norma\phih_{\L\infty W}^2 + \norma\phi_{\L\infty W}^2 \bigr) \, \norma{\phih-\phi}_{\L\infty V}^4
  \non
  \\
  & \leq c \, \norma h_{\L2\Vp}^4\,,
\end{align}
\juerg{where the last inequality follows from} \eqref{stability} and \eqref{contdep}.
At this point, we \pier{come} back to \eqref{persecondaF}\an{,} integrate it with respect to time, choose $\delta$ small enough,
and apply Gronwall's lemma.
By also accounting for \eqref{primaF}, we conclude~that
\Beq
  \norma\psih_{\L\infty H} + \norma\omegah_{\L2V}
  \leq c \, \norma h_{\L2\Vp}^2 \,.
  \label{secondaF}
\Eeq

\step
Conclusion

By revisiting the expressions of $\Lambda_1$ and $\Lambda_2$ and their treatement in the above calculations,
and accounting for the estimates on $\psih$ and~$\omegah$ we have obtained,
we deduce~that
\Beq
  \norma{\Lambda_1}_{\L2H} \leq c \, \norma h_{\L2\Vp}^2
  \aand
  \norma{\Lambda_2}_{\L2H} \leq c \, \norma h_{\L2\Vp}^2 \,.
  \non
\Eeq
Hence, by comparison, first in \eqref{secondah} and then in \eqref{primah},
we deduce~that
\Beq
  \norma\etah_{\L2H}\leq c \, \norma h_{\L2\Vp}^2
  \aand
  \norma{\dt\psih}_{\L2\Wp} \leq c \, \norma h_{\L2\Vp}^2 \,.
  \non
\Eeq
The last inequality, when combined with \eqref{daprimaF} and \eqref{secondaF}, yields \eqref{tesi}, and the proof is complete.
\Edim

\Brem
\label{Moregencost}
The inequality \eqref{tesi} \juerg{also} shows the \Frechet\ differentiability
of the map $u\mapsto w$ as a mapping from $\UR$ into \pier{$\L2V$}, 
where $w$ is the third component of the solution to the state system corresponding to~$u$.
Thus, we could have \an{considered in} the cost functional \eqref{cost} an integral depending on~$w$ of the~type
$\pier{(\alpha_4/2)}\intQ |w - \wQ|^2$ with some target function $\wQ\in\LQ2$ \an{and a nonnegative constant $\alpha_4$}.
The whole theory could have been developed with minor changes, indeed.
\juerg{However,} \an{such a target term is not that relevant from the modeling viewpoint\juerg{. Therefore, its inclusion would be rather a mathematical \pier{extension}, wherefore we decided to not include it,} being aware that it may be handled from a mathematical standpoint.}
\Erem


\subsection{Necessary conditions for optimality}
\label{NC}

By virtue of the \Frechet\ differentiability result given by Theorem~\ref{Frechet},
\an{along} with the continuous embedding $\H1\Wp\cap\L2W\emb\C0H$
and the quadratic structure of the cost functional~\eqref{cost},
we can apply the chain rule to the \an{composite}~mapping
\begin{align}
  & \UR \to (\UR \times \calY) \to \erre
  \quad \hbox{given by} \quad  
  u \mapsto (u,\phi) := (u,\calS(u)) \, \mapsto\,  \calJ(u,\phi) \,.
  \non
\end{align}
Since the control \an{problem} consists in minimizing this mapping on $\Uad$ and $\Uad$ is convex, 
\an{a} necessary condition for some $\ustar$ to be an optimal \an{control} is that
the variational inequality
\Beq
  \alpha_1 \intQ (\phistar-\phiQ) \psi
  + \alpha_2 \iO (\phistar(T)-\phiO) \psi(T)
  + \alpha_3 \intQ \ustar (u-\ustar)
  \geq 0
  \label{badNC}
\Eeq
holds true for every $u\in\Uad$,
where $\psi$ is the first component of the solution $\soluzl$ to the linearized problem
associated with $h:=u-\ustar$ and the components $\phistar$ and $\wstar$ of the state $\soluzstar$ corresponding to~$\ustar$.
However, this condition is not satisfactory\juerg{, since it requires to solve} the linearized system infinitely many times
as $u$ \an{varies} in~$\Uad$.
\an{Let us briefly mention that proceeding as suggested in Remark \ref{Moregencost} would introduce an additional 
term of the form $ \alpha_4\intQ (\wstar- w_Q)\omega $ in \eqref{badNC},
\juerg{with $\omega$ being} the third component of the solution $\soluzl$ to the linearized problem
associated with $h:=u-\ustar$.}

\juerg{As usual, this unpleasant situation} is overcome by the introduction of a proper adjoint problem
associated with the above optimal control and the corresponding state.
This problem consists in looking for a triplet $\soluza$ satisfying
\begin{align}
  & p \in \H1\Wp \cap \L\infty V \cap \L2{\an{{}\Hx4\cap {}}W} , \quad
  q \in \L2W ,
  \non
  \\
  & \aand
  r \in \L2H ,
  \label{regsoluza}
  \\
  & - \< \dt p , v >_W
  - \iO r \, \Delta v
  + \sigma \iO p v
  + \iO \xistar q v
  \pier{{}+{}} \iO \zetastar r v
  \non
  \\
  & = \iO \rho_1 v
  \quad \hbox{\aet, for every $v\in W$,}
  \label{primaa}
  \\
  \separa
  & q + \Delta p = 0
  \quad \aeQ ,
  \label{secondaa}
  \\
  \separa
  & \iO r v
  - \iO \nabla q \cdot \nabla v
  - \nu \iO q v
  - \iO \zetastar q v
  \non
  \\
  & = 0
  \quad \hbox{\aet, for every $v\in V$,}
  \label{terzaa}
  \\
  \separa
  & p(T) = \rho_2,
  \label{cauchya}
\end{align}
\Accorpa\Pbla primaa cauchya
where, for brevity, we have set
\begin{align}
  & \xistar := f''(\phistar) \, \wstar , \quad
  \zetastar := f'(\phistar) , \quad
  \rho_1 := \alpha_1 (\phistar-\phiQ) 
  \non
  \\
  & \aand
  \rho_2 := \alpha_2 (\phistar(T)-\phiO) \,.
  \label{defrho}
\end{align}

\Bthm
\label{Wellposednessa}
\juerg{In addition to the current assumptions, suppose $\phiO\in V$. Then}
the adjoint problem has a unique solution $\soluza$ satisfying~\eqref{regsoluza}.
\Ethm

\Bdim
We first perform a basic estimate and prove uniqueness.
Then we sketch how the same estimate can be used to prove existence.

\step
The basic estimate

We test \eqref{primaa} first by $p$, and then by $q$, and respectively obtain~that
\begin{align}
  & - \< \dt p , p >_W
  - \iO r \Delta p
  + \sigma \iO |p|^2
   = - \iO \xistar q p
  \pier{{}-{}} \iO \zetastar r p
  + \iO \rho_1 p \,,
  \non
  \\
  & - \< \dt p , q >_W
  - \iO r \, \Delta q
  = - \sigma \iO pq
  - \iO \xistar |q|^2
 \pier{{}-{}}  \iO \zetastar r q
  + \iO \rho_1 q\an{\,.}
  \non
\end{align}
Now, we notice that \eqref{regsoluza} and \eqref{secondaa} yield that $\Delta p\in\L2W$.
Thus, we can take the duality beween $\dt p$ and equation \eqref{secondaa} \an{leading to}
\Beq
  \< \dt p , q >_W
  + \< \dt p , \Delta p >_W
  = 0 \,.
  \non
\Eeq
Next, \an{for} a positive constant $K$ whose value is chosen later on,
we test \eqref{secondaa} by $Kq$ and \pier{integrate by parts to infer}~that
\Beq
  K \iO |q|^2 
  = \pier{{} K \iO \nabla q \cdot \nabla p {}}\,.
  \non
\Eeq
Now, we test \eqref{secondaa} by~$r$.
Moreover, we perform an integration by parts in \eqref{terzaa},
write the resulting equality as a differential equation,
multiply by $r$ and integrate over~$\Omega$.
Finally, we test \eqref{terzaa} as it is by~$-q$.
We \juerg{obtain the following identities:}
\begin{align}
  & \iO q r
  + \iO r \Delta p
  = 0
  \non
  \\
  & \iO |r|^2
  + \iO \Delta q \, r
  = \nu \iO q r
  + \iO \zetastar q r
  \non
  \\
  & - \iO r q
  + \iO |\nabla q|^2
  = - \nu \iO \an{|q|^2}
  - \iO \zetastar |q|^2 \,.
  \non
\end{align}
At this point, we add all of the equalities we have \last{listed above} to each other 
and notice several cancellations leading to \juerg{the identity}
\begin{align}
  & - \< \dt p , p >_W
  + \< \dt p , \Delta p >_W
  + \sigma \iO |p|^2
  + K \iO |q|^2
  + \iO |r|^2
  + \iO |\nabla q|^2
  \non
  \\
  & = - \iO \xistar p q
  + \iO \zetastar p r
  + \iO \rho_1 p
	\juerg{{}-\sigma \iO pq
	-\iO \xistar |q|^2
	+\iO \zetastar  q r  {}}
  + \iO \rho_1 q
  \non
  \\
  & \quad {}
  \juerg{{}+ K \iO \nabla q \cdot \nabla p{}}
  + \nu \iO q r 
  + \iO \zetastar q r
  - \nu \iO |q|^2
  - \iO \zetastar |q|^2 \,.
  \label{perbasic}
\end{align}
In order to treat the \lhs, we \pier{first note~that}
\Beq
  - \< \dt p , p >_W
  = - \frac 12 \, \ddt \, \norma p^2 \,.
  \non
\Eeq
On the other hand, 
\an{as a direct consequence of Remark \ref{Bravogiannibis},}
we have~that \an{the second term on the \lhs\ can be interpreted as}
\Beq
  \< \dt p , \Delta p >_W
  = - \frac 12 \, \ddt \, \norma{\nabla p}^2 \,.
  \non
\Eeq
Therefore, by integrating \eqref{perbasic} over $(t,T)$ and observing that
the assumption $\phiO\in V$ and the regularity of $\phistar$ ensure that $\rho_2\in V$,
the sum of the first two terms is given~by
\Beq
  \an{\frac 12 \, \norma{p(t)}^2_V}
  \an{- \frac 12 \, \norma{\rho_2}^2_V\,.}
  \non
\Eeq
As for the time integral of the \rhs, we \pier{recall~\eqref{stability} but} do not detail how to estimate it.
However, it is easy to see that one can play with the Young inequality and choose $K$ large enough
in order to apply the Gronwall lemma. \juerg{In this connection, the terms involving $\zetastar$ are
easy to handle, since $\zetastar$ is bounded. \pier{The} terms involving $\xistar$ require some 
attention. However, we know that $\xistar\in L^\infty(0,T;H)$, and thus we can estimate the only critical term
by means of the compactness inequality \eqref{compact} as follows:}
\begin{align*}
&\juerg{{}-\int_t^T \!\! \iO \xistar |q|^2 \,\le\,\int_t^T\|\xistar(s)\|\,\|q(s)\|_4^2\,ds}\\
&\juerg{{}\le\,c\int_t^T\|q(s)\|_4^2\,ds\,\le \,\frac 14\int_t^T\!\!\iO |\nabla q|^2\,+\,c\,\int_t^T\!\!\iO |q|^2,}
\end{align*}
\juerg{where the last term can be absorbed by choosing a sufficiently large $K$.} 
\an{Thus,} we can conclude~that
\Beq
  \norma p_{\L\infty V}
  + \norma q_{\L2V}
  + \norma r_{\L2H}
  \leq c \, \bigl( \norma{\rho_1}_{\L2H} + \normaV{\rho_2} \bigr) \,.
  \label{basic}
\Eeq

\step
Uniqueness

Recalling that \eqref{basic} has been rigorously obtained under the regularity assumptions \eqref{regsoluza},
and applying it with $\rho_1$ and $\rho_2$ replaced by zero,
we conclude that the solution is $(0,0,0)$.
Since the problem is linear, this proves uniqueness.

\step
Existence

As in the case of the original system \Pbl\ and of the linearized system \Pbll,
one can construct a solution to \Pbla\ by starting from a Faedo\an{--}Galerkin scheme
with the same set of eigenfunctions as considered before.
Then, estimate \eqref{basic} can be \an{rigorously} performed at the discrete level.
\juerg{In fact}, the argument we have used can even be simplified, since the discrete solution is smooth.
Once the analogue of \eqref{basic} is obtained, 
one can derive further estimates by comparison in the discrete equations
(see, e.g., the argument used to prove \eqref{sestastima}).
The formal analogues on our system are 
\begin{align}
  & \norma{\Delta q}_{\L2H} \leq c\,,
  \quad \hbox{whence} \quad
  \norma q_{\L2W} \leq c\,,
  \non
  \\
  & \norma{\Delta p}_{\L2W} \leq c\,,
  \quad \hbox{whence} \quad
  \norma p_{\L2{\Hx4}} \leq c\,,
  \non
  \\
  & \norma{\dt p}_{\L2\Wp} \leq c\,,
  \non
\end{align}
\juerg{which are} obtained by comparison in \juerg{\eqref{terzaa}, \eqref{secondaa}, and \eqref{primaa}}, respectively.
Once estimates like these are proved at the discrete setting,
one can easily let the \last{discretization parameter} tend to infinity and obtain a solution to problem \Pbla.
\Edim

We conclude the paper by proving a \an{first-order} necessary condition for optimality.

\Bthm
\label{GoodNC}
Let $\ustar\in\Uad$ be an optimal control, and let $\soluzstar$ be the corresponding state.
Then
\Beq
  \intQ (\alpha_3 \, \ustar + p) (u-\ustar) \geq 0 \quad \hbox{for every $u\in\Uad$,}
  \label{goodNC}
\Eeq
where $p$ is the first component of the solution $\soluza$ to the corresponding adjoint problem.
In particular, if $\alpha_3>0$, \juerg{then} the optimal control $\ustar$ is the $H$-projection of $-p/\alpha_3$ onto~$\Uad$.
\Ethm

\Bdim
We fix $u$ in $\Uad$ and write the linearized system \accorpa{primal}{cauchyl} 
associated with the optimal control $\ustar$, the corresponding state $\soluzstar$, and the variation $h:=u-\ustar$,
and test the three equations by $p$, $q$, and~$r$, respectively.
More precisely, regarding \eqref{terzal}, on account of the regularity of $\psi$ given in \eqref{regsoluzl},
we replace the first integral by $-\iO\Delta\psi\,v$ and then test by~$r$.
At the same time, we test \eqref{primaa} by $-\psi$, multiply \eqref{secondaa} by $\eta$, and integrate over~$\Omega$,
and \an{test} \eqref{terzaa} by~$\omega$.
Hence, by recalling the definitions of $\xistar$ and $\zetastar$ given in~\eqref{defrho},
we obtain \juerg{a.e. in $(0,T)$} the \an{identities}
\begin{align}
  & \< \dt\psi , p >_W 
  \an{
  - \iO \eta \, \Delta p
  + \sigma \iO \psi p 
  }
  = \iO (u-\ustar) p\,,
  \non
  \\
  \separa
  & \iO \nabla\omega \cdot \nabla q
  + \iO f''(\phistar) \, \wstar  \psi  q
  + \iO f'(\phistar) \, \omega  q
  + \nu \iO \omega q
  = \iO \eta  q\,,
  \non
  \\
  \separa
  & - \iO \Delta\psi \, r + \iO f'(\phistar) \, \psi r
  = \iO \omega r\,,
  \non
  \\
  \separa
  & \< \dt p , \psi >_W
  + \iO r \, \Delta\psi
  - \sigma \iO p \psi
  - \iO f''(\phistar) \, \wstar q \psi
  - \iO f'(\phistar) \, r \psi
  \non
  \\
  & = - \iO \rho_1 \psi\,,
  \non
  \\
  \separa
  & \iO q \eta + \iO \Delta p \, \eta = 0\,,
  \non
  \\
  \separa
  & \iO r \omega
  - \iO \nabla q \cdot \nabla \omega
  - \nu \iO q \omega
  - \iO f'(\phistar) \, q \omega
  = 0 \,.
  \non
\end{align}
At this point, we add all of them to each other and integrate over~$(0,T)$.
Due to obvious cancellations, what remains is~the identity
\Beq	
  \ioT \bigl( \< \dt\psi(t) , p(t) >_W + \< \dt p(t) , \psi(t) >_W \bigr) \, dt
  =  \intQ (u-\ustar) p
  - \intQ \rho_1 \psi \,,
  \non
\Eeq
and an application of the well-known integration-by-parts formula in the framework of the Hilbert triplet $(W,H,\Wp)$, 
combined with \eqref{cauchyl} and \eqref{cauchya}, yields~that
\Beq
  \intQ \rho_1 \psi
  + \iO \rho_2 \psi(T)
  = \intQ (u-\ustar) p \,.
  \non
\Eeq
Therefore, by recalling \eqref{defrho} and using the above equality in \eqref{badNC},
we obtain~\eqref{goodNC}.
\Edim

\juerg{
\Brem
Since the optimal control problem {\bf (P)} is nonconvex, it will usually have many local minima. In this 
connection, recall that a control $u^*\in\Uad$ is called {\em locally optimal for {\bf (P)} in the sense of $L^p$ 
for $p\in[1,+\infty]$} if and
only if there is some $\tau>0$ such that ${\cal J}(u^*,{\cal S}(u^*))\,\le\,{\cal J}(u,{\cal S}(u))$ for all $u\in\Uad$ 
with $\|u-u^*\|_p\le\tau$. It is easily seen that any locally optimal control in the sense of $L^p$ for some $p\in [1,+\infty)$ is
also locally optimal in the sense of $L^\infty$. We now claim that the variational inequality \eqref{goodNC} is valid also
for every control which is locally optimal in the sense of any $p\in [1,+\infty]$. Indeed, \last{it is easily observed} that for 	
every locally optimal control $u^*$ in the sense of $L^\infty$ the variational inequality \eqref{badNC} must be satisfied. 
By the same argument as above, then also \eqref{goodNC} must be valid.  
\Erem
}


\vskip 6mm
\noindent{\bf Acknowledgements}

\noindent
This research was supported by the Italian Ministry of Education, 
University and Research (MIUR): Dipartimenti di Eccellenza Program (2018--2022) 
-- Dept.~of Mathematics ``F.~Casorati'', University of Pavia. 
In addition, {PC and AS gratefully acknowledge some other support 
from the MIUR-PRIN Grant 2020F3NCPX ``Mathematics for industry 4.0 (Math4I4)'' and}
their affiliation to the GNAMPA (Gruppo Nazionale per l'Analisi Matematica, 
la Probabilit\`a e le loro Applicazioni) of INdAM (Isti\-tuto 
Nazionale di Alta Matematica). 
\last{AS has been supported by ``MUR GRANT Dipartimento di Eccellenza'' 2023-2027.}


\End{document}
